\documentclass[12pt]{amsart}
\usepackage{amsmath,amsthm,amsfonts,amscd,amssymb,eucal,latexsym,mathrsfs, appendix, yhmath}
\usepackage[numbers,sort&compress]{natbib}
\usepackage[all,cmtip]{xy}

\newtheorem{theorem}{Theorem}[section]

\newtheorem{lemma}[theorem]{Lemma}
\newtheorem{proposition}[theorem]{Proposition}

\theoremstyle{definition}
\newtheorem{definition}[theorem]{Definition}
\newtheorem{remark}[theorem]{Remark}

\newtheorem{question}[theorem]{Question}

\newcommand{\Zb}{{\mathbb Z}}

\newcommand{\Rb}{{\mathbb R}}
\newcommand{\Nb}{{\mathbb N}}

\newcommand{\cB}{{\mathcal B}}
\newcommand{\cU}{{\mathcal U}}

\newcommand{\cV}{{\mathcal V}}

\newcommand{\cW}{{\mathcal W}}

\newcommand{\cE}{{\mathcal E}}
\newcommand{\cA}{{\mathcal A}}
\newcommand{\cD}{{\mathcal D}}
\newcommand{\cZ}{{\mathcal Z}}

\newcommand{\diam}{{\rm diam}}
\newcommand{\sym}{{\rm Sym}}
\newcommand{\map}{{\rm Map}}

\allowdisplaybreaks
\numberwithin{equation}{section}

\begin{document}
\title[Topological pressure and the variational principle]{Topological pressure and the variational principle for actions of sofic groups}
\author{Nhan-Phu Chung}
\address{\hskip-\parindent
Nhan-Phu Chung, Department of Mathematics, SUNY at Buffalo,
Buffalo NY 14260-2900, U.S.A.}
\email{phuchung@buffalo.edu}


\date{April 07, 2012}
\begin{abstract} We introduce topological pressure for continuous actions of countable sofic groups on compact metrizable spaces. This generalizes the classical topological pressure for continuous actions of countable amenable groups on such spaces. We also establish the variational principle for topological pressure in this sofic context.
\end{abstract}
\maketitle
\section{Introduction}
Starting from ideas in the statistical mechanics of lattice systems, in \cite{Ruelle73} Ruelle introduced topological pressure of a continuous function for actions of the groups $\Zb^n$ on compact spaces and established the variational principle for topological pressure in this context when the action is expansive and satisfies the specification condition. Later, Walters \cite{Walters75} dropped these assumptions when he proved the variational principle for a $\Zb^+$-action. A shorter and elegant proof of the variational principle for $\Zb^n_+$-actions was given by Misiurewicz \cite{Misiurewicz76}. Stepin and Tagi-Zade \cite{StepinTagi80}, Moullin Ollagnier and Pinchon \cite{MouPinchon82, Mou85}, Tempelman \cite{Tempelman84, Tempelman92} extended the variational principle to the case when $\Zb^n$ is replaced by any countable amenable group.

From a viewpoint of dimension theory, Pesin and Pitskel' \cite{PesinPitskel84} introduced another way to define topological pressure for continuous functions on noncompact sets in the case of $\Zb$-actions. For more information and references in this direction, see \cite{Pesin97}.

The notion of a sofic group was first introduced by Gromov \cite{Gromov99}. All countable amenable groups and residually finite groups are sofic. It is unknown whether every countable group is sofic. We refer readers to \cite{CeccheriniCoornaert10, ElekSzabo05, ElekSzabo06, Pestov08, Thom08, Weiss00} for details on sofic groups.

In 2008, in a remarkable result, Lewis Bowen \cite{Bowen10} defined sofic entropy for measure-preserving actions of countable sofic groups on standard probability measure spaces admitting a generating partition with finite entropy. Recently, in \cite{KerrLi11a, KerrLi10}, via an operator algebraic method, David Kerr and Hanfeng Li extended Bowen's sofic measure entropy to all measure-preserving actions of countable sofic groups on standard probability measure spaces, and defined sofic topological entropy for continuous actions of countable sofic groups on compact metrizable spaces. They also established the variational principle between sofic measure entropy and sofic topological entropy \cite{KerrLi11a}. In the case of amenable groups, the sofic entropies coincide with the classical entropies \cite{Bowen10a, KerrLi10}. After that, the approach of Kerr-Li \cite{KerrLi11a, KerrLi10} for continuous actions of countable sofic groups on compact metrizable spaces has been applied to study mean dimension \cite{Li11} and local entropy theory \cite{Zhang11} in the sofic context.

Given Kerr-Li's work, it is natural to ask how to define topological pressure of a continuous function for actions of countable sofic groups on compact metrizable spaces and whether it coincides with the classical topological pressure for actions of countable amenable groups on such spaces. Furthermore, one might ask whether there exists a relation between sofic topological pressure and sofic measure entropy via a variational principle.

The goal of this paper is to answer all of these questions. We organize this paper as follows. We define the sofic topological pressure $P_{\Sigma}(f,X,G)$ and establish some basic properties of it in Section 2. In Section 3, we recall the definition of classical topological pressure $P(f,X,G)$ for actions of countable amenable groups and prove our first main result:
\begin{theorem}
\label{T-sofic pressure is equal to amenable pressure}
Let $G$ be a countable amenable group acting continuously on a compact metrizable space $X$. Let $\Sigma$ be a sofic approximation sequence for $G$ and $f$ be a real valued continuous function on $X$. Then $P_{\Sigma}(f,X,G)= P(f,X,G)$.
\end{theorem}
In Section 4, we recall the definition of sofic measure entropy $h_{\Sigma,\mu}(X,G)$ and prove our second main result about the variational principle for sofic topological pressure. The variational principle for topological pressure is well known when the acting group $G$ is amenable. For example, see \cite[Theorem 9.10]{Walters82} for the case $G=\Zb$ and \cite[Theorem 5.2.7]{Mou85} for the case $G$ is a countable amenable group.
\begin{theorem}
\label{T-variational principle for topological pressure}
Let $\alpha$ be a continuous action of a countable sofic group $G$ on a compact metrizable space $X$. Let $\Sigma$ be a sofic approximation sequence for $G$ and $f$ be a real valued continuous function on $X$. Then $$P_{\Sigma}(f,X,G)=\sup\Big\{h_{\Sigma,\mu}(X,G)+\int_X fd\mu:\mu\in M_G(X)\Big\},$$
where $M_G(X)$ is the set of $G$-invariant Borel probability measures on $X$. In particular, if $P_{\Sigma}(f,X,G)\neq -\infty$ then $M_G(X)$ is nonempty.
\end{theorem}
To illustrate an example, we compute the sofic topological pressure and find some equilibrium state for some function on Bernoulli shifts in Section 5. Finally, in Section 6, we describe some properties of topological pressure and give a sufficient condition for a finite signed measure to be a member of $M_G(X)$, using topological pressure.

To finish the Introduction, we recall the definitions of sofic groups, separated sets, and spanning sets and fix some notations.

For each $d\in\Nb$, we denote by $[d]$ the set $\{1,...,d\}$ and $\sym(d)$ the permutation group of $[d]$.

For every real number $y$, we denote by $\lfloor y\rfloor$ the largest integer which is less than or equal to $y$.

Let $G$ be a countable group. We say that $G$ is \textit{sofic} if there is a sequence $\Sigma=\{\sigma_i: G \rightarrow \sym(d_i), d_i\in \Nb\}_{i\in \Nb}$ such that
\begin{enumerate}
\item $\lim_{i\rightarrow \infty}\frac{1}{d_i}|\{a\in [d_i]:\sigma_{i,s}\sigma_{i,t}(a)=\sigma_{i,st}(a)\}|=1$ for all $s,t\in G$,

\item $\lim_{i\rightarrow \infty}\frac{1}{d_i}|\{a\in[d_i]:\sigma_{i,s}(a)\neq \sigma_{i,t}(a)\}|=1$ for all distinct $s,t\in G$,

\item $\lim_{i\rightarrow \infty}d_i=\infty$.
\end{enumerate}
Such a sequence is called a \textit{sofic approximation sequence} for $G$.
Note that when $G$ is infinite, the condition (3) is a consequence of the condition (2).

Let $(Y,\rho)$ be a pseudometric space and $\varepsilon>0$. A subset $A$ of $Y$ is called $(\rho,\varepsilon)$-\textit{separated} if $\rho(x,y)\geq \varepsilon$ for all distinct $x,y\in A$, and $(\rho,\varepsilon)$-\textit{spanning} if for every $y\in Y$ we can find an $x\in A$ such that $\rho(x,y)<\varepsilon$. We denote by $N_\varepsilon(Y,\rho)$ the maximal cardinality of a finite $(\rho,\varepsilon)$-separated subset of $Y$.

 Throughout this paper, the space $X$ is always compact metrizable and $G$ is always a countable sofic group with the identity element $e$. We denote by $C(X)$ the set of all real valued continuous functions on $X$. A continuous action $\alpha$ of $G$ on a compact metrizable space $X$ induces an action of $G$ on $C(X)$ as follows: for $g\in C(X)$ and $s\in G$, the function $\alpha_s(g)$ is given by $x\mapsto g(s^{-1}x)$. Given a map $\sigma:G\rightarrow \sym(d)$ for some $d\in\Nb$, for $s\in G, x\in X,$ and $a\in [d]$ we will for convenience denote $\alpha_s(x)$ and $\sigma_s(a)$ by $sx$ and $sa$ respectively.

Let $\rho$ be a continuous pseudometric on $X$. For any $d\in\Nb$, we define the pseudometrics $\rho_2,\rho_\infty$ on the set of all maps from $[d]$ to $X$ as follows:
\begin{align*}
\rho_2(\psi,\varphi) & =\Big(\frac{1}{d}\sum_{i=1}^d(\rho(\psi(i),\varphi(i)))^2\Big)^{1/2},\\
\mbox{~~~ and ~~~~~}\rho_\infty(\psi,\varphi) & = \max_{1\leq i\leq d}\rho(\psi(i),\varphi(i)).
\end{align*}
For every subset $J$ of $[d]$, we define on the set of maps from $[d]$ to $X$ the pseudometric
$$
\rho_{J,\infty} (\psi , \varphi ): = \rho_{\infty} (\psi|_J,  \varphi|_J).
$$

\noindent{\bf Acknowledgements:} I am grateful to my advisor, Prof. Hanfeng Li for introducing me to the subject and his continuous guidance, support and encouragement. I also thank Guohua Zhang and the referee for helpful comments. 
\section{Sofic Topological Pressure}
In this section, we will define the topological pressure of a continuous function for actions of countable sofic groups on compact metrizable spaces and establish some basic properties of it.

Let $\alpha$ be a continuous action of a countable sofic group $G$ on a compact metrizable space $X$. Let $f$ be a real valued continuous function on $X$, $\rho$ a continuous pseudometric on $X$ and $\Sigma$ a sofic approximation sequence of $G$. Let $F$ be a nonempty finite subset of $G$ and $\delta>0$. Let $\sigma$ be a map from $G$ to $\sym(d)$ for some $d\in \Nb$.
Now we recall the definition of $\map(\rho,F,\delta,\sigma)$.
\begin{definition}
 We define $\map(\rho,F,\delta,\sigma)$ to be the set of all maps $\varphi:[d]\rightarrow X$ such that $\max_{s\in F}\rho_2(\alpha_s\circ \varphi, \varphi\circ\sigma_s)<\delta$.
\end{definition}

The space $\map(\rho,F,\delta,\sigma)$ appeared first in \cite[Section 2]{KerrLi10}, and has been applied to study sofic entropies \cite{KerrLi11a}, sofic mean dimension \cite{Li11}, local entropy theory \cite{Zhang11}.

\begin{definition}
Let $\varepsilon>0$. We define $$M^\varepsilon_{\Sigma,\infty}(f,X,G,\rho,F,\delta,\sigma)=\sup_{\cE}\sum_{\varphi\in \cE} \exp\Big(\sum_{a=1}^df(\varphi(a))\Big),$$
where $\cE$ runs over $(\rho_\infty,\varepsilon)$-separated subsets of $\map(\rho,F,\delta,\sigma)$. Of course, the value of the right hand side doesn't change if $\cE$ runs over maximal $(\rho_\infty,\varepsilon)$-separated subsets of $\map(\rho,F,\delta,\sigma)$.
\end{definition}
Now we define the sofic topological pressure of $f$.
\begin{definition}
We define
\begin{eqnarray*}
P^\varepsilon_{\Sigma,\infty}(f,X,G,\rho,F,\delta)&=&\limsup_{i\rightarrow \infty}\frac{1}{d_i}\log M^\varepsilon_{\Sigma,\infty}(f,X,G,\rho,F,\delta,\sigma_i),\\ P^\varepsilon_{\Sigma,\infty}(f,X,G,\rho,F)&=&\inf_{\delta>0}P^\varepsilon_{\Sigma,\infty}(f,X,G,\rho,F,\delta),\\
P^\varepsilon_{\Sigma,\infty}(f,X,G,\rho)&=&\inf_F P^\varepsilon_{\Sigma,\infty}(f,X,G,\rho,F),\\
P_{\Sigma,\infty}(f,X,G,\rho)&=&\sup_{\varepsilon>0}P^\varepsilon_{\Sigma,\infty}(f,X,G,\rho),
\end{eqnarray*}
where $F$ in the third line runs over the nonempty finite subsets of $G$.

 If $\map(\rho,F,\delta,\sigma_i)=\emptyset$ for all large enough $i$, we set $P^\varepsilon_{\Sigma,\infty}(f,X,G,\rho,F,\delta)=-\infty$.

Similarly, we define $M^\varepsilon_{\Sigma,2}(f,X,G,\rho,F,\delta,\sigma_i), P^\varepsilon_{\Sigma,2}(f,X,G,\rho,F,\delta), P^\varepsilon_{\Sigma,2}(f,X,G,\rho,F),$\\ $P^\varepsilon_{\Sigma,2}(f,X,G,\rho)$ and $P_{\Sigma,2}(f,X,G,\rho)$ using $\rho_2$ in place of $\rho_\infty$.
 \begin{remark}
 \label{R-topological entropy is a special case of topological pressure}
When $f=0$, $P_{\Sigma,\infty}(0,X,G,\rho)$ is the sofic topological entropy $h_{\Sigma,\infty}(X,G,\rho)$, as defined in \cite[Section 2]{KerrLi10} and originating in another equivalent form in \cite[Section 4]{KerrLi11a}.
 \end{remark}

 Now we prove that the definition of sofic topological pressure does not depend on the choice of $\rho_2$ and $\rho_\infty$.
 \begin{lemma}
Let $\rho$ be a continuous pseudometric on $X$ such that $f$ is continuous with respect to $\rho$. Then $$P_{\Sigma,2}(f,X,G,\rho)=P_{\Sigma,\infty}(f,X,G,\rho).$$
 \end{lemma}
\begin{proof}
Since $\rho_\infty\geq\rho_2$, $P_{\Sigma,2}(f,X,G,\rho)\leq P_{\Sigma,\infty}(f,X,G,\rho).$

Now we prove $P_{\Sigma,\infty}(f,X,G,\rho)\leq P_{\Sigma,2}(f,X,G,\rho).$

Let $\theta>0$. Let $\varepsilon'>0$ be such that $|f(x)-f(y)|<\theta$ whenever $x,y\in X$ with $\rho(x,y)<\sqrt{\varepsilon'}$.
Let $\varepsilon>0$, which we will determine later.
It suffices to prove that
$$P^{2\sqrt{\varepsilon'}}_{\Sigma,\infty}(f,X,G,\rho,F,\delta)\leq P^\varepsilon_{\Sigma,2}(f,X,G,\rho,F,\delta)+4\theta,$$
for any $\delta>0$ and nonempty finite subset $F$ of $G$. Let $\delta>0$, $F$ be a nonempty finite subset of $G$ and $\sigma$ be a map from $G$ to $\sym(d)$ for some $d\in\Nb$.

Let $\cE$ be a $(\rho_\infty,2\sqrt{\varepsilon'})$-separated subset of $\map(\rho,F,\delta,\sigma)$ such that $$M_{\Sigma,\infty}^{2\sqrt{\varepsilon'}}(f,X,G,\rho, F,\delta,\sigma)\leq 2\cdot\sum_{\varphi\in \cE} \exp\Big(\sum_{i=1}^df(\varphi(i))\Big).$$

Let $\cB$ be a maximal $(\rho_2,\varepsilon)$-separated subset of $\cE$. Then $\cE=\bigcup_{\varphi\in \cB}(\cE\cap B_\varphi)$, where $B_\varphi=\{\psi\in X^{[d]}: \rho_2(\varphi,\psi)<\varepsilon\}.$

Let $\varphi\in \cB$. Let us estimate how many elements are in $\cE\cap B_\varphi$.
Let $Y_{\varepsilon'}$ be a maximal $(\rho,\sqrt{\varepsilon'})$-separated subset of $X$.

For each $\psi\in \cE\cap B_\varphi$, we denote by $\Lambda_\psi$ the set of all $a\in[d]$ such that $\rho(\varphi(a),\psi(a))<\sqrt{\varepsilon'}$. Then $|\Lambda_\psi|\geq (1-\frac{\varepsilon^2}{\varepsilon'})d$. We enumerate the elements of $\{\Lambda_\psi:\psi\in \cE\cap B_\varphi\}$ as $\Lambda_{\varphi,1},...,\Lambda_{\varphi,\ell_\varphi}$.
Then $\cE\cap B_\varphi=\bigsqcup_{j=1}^{\ell_\varphi} \cV_j$, where $\cV_j=\{\psi\in \cE\cap B_\varphi:\Lambda_\psi=\Lambda_{\varphi,j}\}$, for every $j=1,...,\ell_\varphi$.

For every $j=1,...,\ell_\varphi$, set $\Lambda_{\varphi,j}^c=[d]\setminus \Lambda_{\varphi,j}$. Since $Y_{\varepsilon'}$ is a $(\rho,\sqrt{\varepsilon'})$-spanning subset of $X$, for every $\psi \in \cV_j$, we can find $f_\psi\in Y_{\varepsilon'}^{\Lambda_{\varphi,j}^c}$ such that $\rho_\infty(\psi|_{\Lambda_{\varphi,j}^c}, f_\psi)<\sqrt{\varepsilon'}$. Then there exists $\cA\subset \cV_j$ such that $|\cV_j|\leq |Y_{\varepsilon'}|^{|\Lambda_{\varphi,j}^c|}|\cA|$ and $f_\psi$ is the same, say $f$, for every $\psi\in \cA$. Then $$\rho_\infty(\psi|_{\Lambda_{\varphi,j}^c}, \psi'|_{\Lambda_{\varphi,j}^c})\leq \rho_\infty(\psi|_{\Lambda_{\varphi,j}^c}, f)+\rho_\infty(f, \psi'|_{\Lambda_{\varphi,j}^c})<2\sqrt{\varepsilon'},$$
for any $\psi,\psi'\in \cA$.
Since $\cA$ is a $(\rho_\infty,2\sqrt{\varepsilon'})$-separated set, we get $\psi=\psi'$. Thus $|\cA|\leq 1$, and hence $|\cV_j|\leq |Y_{\varepsilon'}|^{|\Lambda_{\varphi,j}^c|}|\cA|\leq |Y_{\varepsilon'}|^{\frac{\varepsilon^2}{\varepsilon'} d}$.

By Stirling's approximation formula, $\frac{\varepsilon^2}{\varepsilon'} d \binom{d}{\frac{\varepsilon^2}{\varepsilon'} d}$ is less than $\exp(\beta d)$ for some $\beta>0$ depending on $\varepsilon$ but not on $d$ when $d$ is large enough with $\beta\rightarrow 0$ as $\varepsilon\rightarrow 0$. Since $\sum_{j=0}^{\lfloor \frac{\varepsilon^2}{\varepsilon'} d \rfloor} \binom{d}{j}\leq \frac{\varepsilon^2}{\varepsilon'} d \binom{d}{\frac{\varepsilon^2}{\varepsilon'} d}$, when $d$ is large enough we have that the number of subsets of $[d]$ of cardinality at least $(1-\frac{\varepsilon^2}{\varepsilon'})d$ is at most $\exp(\beta d)$. Therefore,
$$|\cE\cap B_\varphi|\leq \ell_\varphi|Y_{\varepsilon'}|^{\frac{\varepsilon^2}{\varepsilon'} d}\leq \exp(\beta d)|Y_{\varepsilon'}|^{\frac{\varepsilon^2}{\varepsilon'} d}.$$

Since $f$ is continuous on $X$, there exists $Q>0$ such that $|f(x)|\leq Q$ for all $x\in X$. Hence
\begin{eqnarray*}
&& M_{\Sigma,\infty}^{2\sqrt{\varepsilon'}}(f,X,G,\rho, F,\delta,\sigma)\\
& \leq&2\cdot\sum_{\varphi\in \cE} \exp\Big(\sum_{i=1}^df(\varphi(i))\Big)\\
& \leq& 2\cdot\sum_{\varphi\in \cB}\sum_{\psi\in \cE\cap B_\varphi} \exp\Big(\sum_{i=1}^df(\psi(i))\Big)\\
& =& 2\cdot\sum_{\varphi\in \cB}\sum_{\psi\in \cE\cap B_\varphi}  \exp\Big(\sum_{i=1}^df(\varphi(i))\Big)\exp\Big(\sum_{i\in \Lambda_\psi}(f(\psi(i))-f(\varphi(i)))\Big)\\
&&\exp\Big(\sum_{i\notin \Lambda_\psi}(f(\psi(i))-f(\varphi(i)))\Big)\\
& \leq& 2\cdot \sum_{\varphi\in \cB}\sum_{\psi\in \cE\cap B_\varphi}\exp\Big(\sum_{i=1}^df(\varphi(i))\Big)\exp(\theta d)\exp(2Q\frac{\varepsilon^2}{\varepsilon'}  d)\\
&\leq& 2\cdot\sum_{\varphi\in \cB}|Y_{\varepsilon'}|^{\frac{\varepsilon^2}{\varepsilon'} d} \exp(\beta d) \exp\Big(\sum_{i=1}^df(\varphi(i))\Big)
\exp(\theta d+2Q\frac{\varepsilon^2}{\varepsilon'}  d)\\
&\leq& 2\cdot|Y_{\varepsilon'}|^{\frac{\varepsilon^2}{\varepsilon'} d} \exp(\beta d+\theta d+2Q\frac{\varepsilon^2}{\varepsilon'}  d)M_{\Sigma,2}^{\varepsilon}(f,X,G,\rho, F,\delta,\sigma). \\
\end{eqnarray*}
 Thus $P_{\Sigma,\infty}^{2\sqrt{\varepsilon'}}(f,X,G,\rho, F,\delta)\leq P_{\Sigma,2}^{\varepsilon}(f,X,G,\rho, F,\delta)+\frac{\varepsilon^2}{\varepsilon'}  \log N_{\sqrt{\varepsilon'}}(X,\rho)+\beta+\theta+2Q\frac{\varepsilon^2}{\varepsilon'}$.
   We choose $\varepsilon>0$ small enough, not depending on $\delta$ and $F$ such that $\beta<\theta, 2Q\frac{\varepsilon^2}{\varepsilon'}<\theta$ and $\frac{\varepsilon^2}{\varepsilon'}\log N_{\sqrt{\varepsilon'}}(X,\rho)<\theta$. Then
 $P_{\Sigma,\infty}^{2\sqrt{\varepsilon'}}(f,X,G,\rho, F,\delta)\leq P_{\Sigma,2}^{\varepsilon}(f,X,G,\rho, F,\delta)+4\theta$, for all $\delta>0$ and nonempty finite subset $F$ of $G$, as desired.
\end{proof}
\end{definition}
A continuous pseudometric $\rho$ on $X$ is called \textit{dynamically generating} if for any distinct points $x,y\in X$, there exists $s\in G$ such that $\rho(sx,sy)>0$.
The following two lemmas will show that the quantity $P_{\Sigma,\infty}(f,X,G,\rho)$ does not depend on the choice of compatible metric $\rho$ and furthermore it also does not depend on the dynamically generating continuous pseudometric of $X$ with respect to which $f$ is continuous. Thus, we shall write the topological pressure of $f$, $P_{\Sigma,\infty}(f,X,G,\rho) \mbox{ (or }P_{\Sigma,2}(f,X,G,\rho)$), where $\rho$ is a compatible metric on $X$ or a dynamically generating continuous pseudometric on $X$ with respect to which $f$ is continuous, as $P_{\Sigma}(f,X,G)$.
\begin{lemma}
\label{L-Topological pressure does not depend on metrics}
Let $\rho$ and $\rho'$ be compatible metrics on $X$. Then  $P_{\Sigma,\infty}(f,X,G,\rho)=P_{\Sigma,\infty}(f,X,G, \rho')$.
\begin{proof}
 Let $\varepsilon>0$. We choose $\varepsilon'>0$ such that for any $x,y\in X$ with $\rho'(x,y)<\varepsilon'$, one has $\rho(x,y)<\varepsilon$. Let $F$ be a nonempty finite subset of $G$ and $\delta>0$. From the proof in Lemma 2.4 of \cite{Li11}, there exists $\delta'>0$ such that for any map $\sigma$ from $G$ to $\sym(d)$ for some $d\in \Nb$ one has $\map(\rho, F,\delta',\sigma)\subset \map(\rho', F,\delta,\sigma)$. Then any $(\rho_\infty,\varepsilon)$-separated subset of $\map(\rho, F,\delta',\sigma)$ is also a $(\rho'_\infty,\varepsilon')$-separated subset of $\map(\rho', F,\delta,\sigma)$. Thus $$P^\varepsilon_{\Sigma,\infty}(f,X,G,\rho,F)\leq P^\varepsilon_{\Sigma,\infty}(f,X,G,\rho,F,\delta')\leq P^{\varepsilon'}_{\Sigma,\infty}(f,X,G,\rho',F,\delta),$$ and hence $P^\varepsilon_{\Sigma,\infty}(f,X,G,\rho,F)\leq P^{\varepsilon'}_{\Sigma,\infty}(f,X,G,\rho',F)$.

So $P_{\Sigma,\infty}(f,X,G,\rho)\leq P_{\Sigma,\infty}(f,X,G,\rho')$.

Similarly, we also have $P_{\Sigma,\infty}(f,X,G,\rho')\leq P_{\Sigma,\infty}(f,X,G,\rho)$.
\end{proof}
\end{lemma}
\begin{lemma}
Let $\rho$ be a dynamically generating continuous pseudometric on $X$ with respect to which $f$ is continuous. Enumerate the elements of $G$ as $s_1=e,s_2,\dots$. Define a new continuous pseudometric $\rho'$ on $X$ by $\rho'(x,y)=\sum_{n=1}^{\infty}\frac{1}{2^n}\rho(s_nx,s_ny)$ for all $x,y\in X$. Then $\rho'$ is a compatible metric on $X$ and $$P_{\Sigma,\infty}(f,X,G,\rho)=P_{\Sigma,\infty}(f,X,G,\rho').$$
\end{lemma}
\begin{proof}
 Since $\rho$ is dynamically generating, $\rho'$ separates the points of $X$. If we denote by $\tau$ the original topology on $X$, and by $\tau'$ the topology on $X$ induced by $\rho'$, then the identity map $Id:(X, \tau)\rightarrow (X, \tau')$ is continuous. Since $(X, \tau')$ is Hausdorff and $(X, \tau)$ is compact, $Id$ is a homeomorphism. Thus $\rho'$ is a compatible metric on $X$.

Let $\varepsilon>0$. Similar to the proof of \cite[Lemma 4.3]{Li11}, one has $P^{\varepsilon}_{\Sigma,\infty}(f,X,G,\rho)\leq P_{\Sigma,\infty}^{\varepsilon/2}(f,X,G,\rho')$. Thus, $P_{\Sigma,\infty}(f,X,G,\rho)\leq P_{\Sigma,\infty}(f,X,G,\rho').$

Now we will prove $P_{\Sigma,\infty}(f,X,G,\rho')\leq P_{\Sigma,\infty}(f,X,G,\rho)$. It suffices to prove that $P_{\Sigma,\infty}(f,X,G,\rho')\leq P_{\Sigma,\infty}(f,X,G,\rho)+3\theta$ for any $\theta>0$. Let $\theta>0$. Let $\varepsilon'>0$ such that $|f(x)-f(y)|<\theta$ whenever $x,y\in X$ with $\rho(x,y)<\varepsilon'$. It suffices to prove that for any $0<\varepsilon<\varepsilon'$,
$$P^{4\varepsilon}_{\Sigma,\infty}(f,X,G,\rho')\leq P^\varepsilon_{\Sigma,\infty}(f,X,G,\rho)+3\theta.$$

Let $0<\varepsilon<\varepsilon'$. Choose $k\in\Nb$ such that $\diam(X,\rho)/2^k<\varepsilon/2$. Let $F$ be a finite subset of $G$ containing $\{s_1,...,s_k\}$. Let $\delta>0$ be small enough which we will determine later.
Put $\delta'=\delta/2$.  It suffices to prove that $$P^{4\varepsilon}_{\Sigma,\infty}(f,X,G,\rho',F,\delta')\leq P^\varepsilon_{\Sigma,\infty}(f,X,G,\rho,F,\delta)+3\theta.$$

Let $\sigma:G \rightarrow \sym(d)$ be a good enough sofic approximation of $G$, for some $d\in\Nb$.
Since $\rho'_2(\varphi,\psi)\geq \frac{1}{2}\rho_2(\varphi,\psi) $ for all maps $\varphi,\psi:[d]\rightarrow X$, we have $\map(\rho',F,\delta',\sigma)\subset \map(\rho,F,\delta,\sigma)$.

Let $\cE$ be a $(\rho'_\infty,4\varepsilon)$-separated subset of $\map(\rho',F,\delta',\sigma)$ such that $$M^{4\varepsilon}_{\Sigma,\infty}(f,X,G,\rho',F,\delta',\sigma)\leq 2\cdot\sum_{\varphi\in \cE}\exp\Big(\sum_{i=1}^df(\varphi(i))\Big).$$

For each $\varphi\in \cE$ we denote by $\Lambda_\varphi$ the set of all $a\in[d]$ such that $$\max_{s\in F}\rho(\varphi(sa),s\varphi(a))<\sqrt{\delta}.$$ Then $|\Lambda_\varphi|\geq (1-|F|\delta)d$.
We enumerate the elements of $\{\Lambda_\varphi:\varphi\in \cE\}$ as $\Lambda_1,...,\Lambda_\ell$.
Then $\cE=\bigsqcup_{j=1}^\ell \cV_j$, where $\cV_j=\{\varphi\in \cE:\Lambda_\varphi=\Lambda_j\}$, for every $j=1,...,\ell$. Let $Y$ be a maximal $(\rho',2\varepsilon)$-separated subset of $X$. Choose $\delta>0$ such that $\sqrt\delta<\varepsilon/4$ and $|Y|^{|F|\delta}<\exp(\theta)$.

\textbf{Claim}: For any $j=1,...,\ell$, and any $\varphi\in \cV_j$, one has $$|\cV_j\cap B_\varphi|\leq |Y|^{|F|\delta d},$$
where $B_\varphi:=\{\psi\in X^{[d]}:\rho_\infty(\varphi,\psi)<\varepsilon\}$.

A proof of this Claim can be found in the proof of \cite[Lemma 4.3]{Li11}.

By Stirling's approximation formula, $|F|\delta d \binom{d}{|F|\delta d}$ is less than $\exp(\beta d)$ for some $\beta>0$ depending on $\delta$ and $|F|$ but not on $d$ when $d$ is large enough with $\beta\rightarrow 0$ as $\delta\rightarrow 0$. Since $\sum_{j=0}^{\lfloor |F|\delta d\rfloor} \binom{d}{j}\leq |F|\delta d \binom{d}{|F|\delta d}$, when $d$ is large enough we have that the number of subsets of $[d]$ of cardinality at least $(1-|F|\delta)d$ is at most $\exp(\beta d)$. Choose $\delta$ such that $\beta<\theta$. Then, when $d$ is large enough, $\ell\leq \exp(\beta d)\leq \exp(\theta d)$.

For each $j=1,...,\ell$, let $\cB_j$ be a maximal $(\rho_\infty,\varepsilon)$-separated subset of $\cV_j$. Then for any $j=1,...,\ell$, one has $\cV_j=\bigcup _{\varphi\in \cB_j} (\cV_j\cap B_\varphi)$. Thus
\begin{eqnarray*}
&&M_{\Sigma,\infty}^{4\varepsilon}(f,X,G,\rho', F,\delta',\sigma)\\
&\leq&2\cdot\sum_{\varphi\in \cE} \exp\Big(\sum_{i=1}^df(\varphi(i))\Big)\\
&=& 2\cdot\sum_{j=1}^\ell\sum_{\varphi\in \cV_j} \exp\Big(\sum_{i=1}^df(\varphi(i))\Big)\\
&\leq& 2\cdot\sum_{j=1}^\ell\sum_{\varphi\in \cB_j}\sum_{\psi\in \cV_j\cap B_\varphi}  \exp\Big(\sum_{i=1}^d(f(\psi(i))-f(\varphi(i)))\Big)\exp\Big(\sum_{i=1}^df(\varphi(i))\Big)\\
&\leq& 2\cdot\sum_{j=1}^\ell\sum_{\varphi\in \cB_j}\sum_{\psi\in \cV_j\cap B_\varphi}  \exp(\theta d)\exp\Big(\sum_{i=1}^df(\varphi(i))\Big)\\
&\leq& 2\cdot\sum_{j=1}^\ell\sum_{\varphi\in \cB_j}|Y|^{|F|\delta d}  \exp(\theta d)\exp\Big(\sum_{i=1}^df(\varphi(i))\Big)\\
&\leq& 2\cdot\sum_{j=1}^\ell|Y|^{|F|\delta d} \exp(\theta d)M_{\Sigma,\infty}^{\varepsilon}(f,X,G,\rho, F,\delta,\sigma) \\
&=& 2\cdot\ell|Y|^{|F|\delta d} \exp(\theta d)M_{\Sigma,\infty}^{\varepsilon}(f,X,G,\rho, F,\delta,\sigma)\\
&\leq& 2\cdot \exp(3\theta d)M_{\Sigma,\infty}^{\varepsilon}(f,X,G,\rho, F,\delta,\sigma).
\end{eqnarray*}
Therefore, $P^{4\varepsilon}_{\Sigma,\infty}(f,X,G,\rho',F,\delta')\leq P^\varepsilon_{\Sigma,\infty}(f,X,G,\rho,F,\delta)+3\theta$.
\end{proof}
\section{Topological pressure in the amenable case}
The purpose of this section is to prove Theorem 1.1.

We begin this section by recalling the classical definition of topological pressure in Section 5 of \cite{Mou85}. A countable group $G$ is said to be \textit{amenable} if there exists a F$\o$lner sequence, which is a sequence $\{ F_i \}_{i=1}^\infty$ of nonempty finite subsets of $G$ such that
$\frac{|sF_i \Delta F_i |}{|F_i|} \to 0$ as $i\to\infty$ for all $s\in G$. We refer the readers to \cite{Paterson88} for details on amenable groups.

Let $G$ be a countable amenable group and $\alpha$ a continuous action of $G$ on a compact metrizable space $X$.
Let $\rho$ be a compatible metric on $X,f\in C(X)$, $F$ a nonempty finite subset of $G$ and $\delta>0$. We define the metric $\rho_F$ on $X$ by $\rho_F(x,y)=\max_{s\in F}\rho(sx,sy)$. An open cover $\cU$ of $X$ is said to be of order $(F,\delta)$ if for any $U\in\cU,$ and $x,y\in U,$ one has $\max_{s\in F}\rho(sx,sy)<\delta$.
We define
$$P_1(F,f,\delta)=\inf_\cU \sum_{U\in\cU} \sup_{x\in U}\exp \big(\sum_{s\in F}f(\alpha_s(x))\big),$$ where $\cU$ runs over the set of all finite open covers of order $(F,\delta)$. By the Ornstein-Weiss lemma in Theorem 6.1 of \cite{LindWeiss}, for any $\delta>0$ the quantities $$\frac{1}{|F|}\log P_1(F,f,\delta)$$ converge to a number, denoted by $p_1(f,\delta)$, as $F$ becomes more and more left invariant in the sense
that for every $\varepsilon > 0$ there are a nonempty finite set $K\subseteq G$ and a $\delta' > 0$
such that $$\big|\frac{1}{|F|}\log P_1(F,f,\delta)-p_1(f,\delta)\big|<\varepsilon,$$ for any nonempty finite subset $F$ of $G$ satisfying $|KF\Delta F| \leq \delta' |F|$.
The topological pressure of $f$ is defined as $\sup_{\delta>0}p_1(f,\delta)$ and does not depend on the choice of compatible metric $\rho$. We denote the topological pressure of $f$ by $P(f,X,G)$.

For any nonempty finite subset $F$ of $G,\varepsilon>0$ and any compatible metric $\rho$ on $X$, define
$$K_\varepsilon(f,X,G,\rho,F)=\sup_{\cD}\sum_{x\in \cD} \exp\big(\sum_{s\in F}f(\alpha_s(x))\big),$$ where $\cD$ runs over $(\rho_{F},\varepsilon)$-separated subsets of $X$.
Given a F{\o}lner sequence $\{F_n\}_{n=1}^\infty$ of $G$, the topological pressure of $f$ can be alternatively expressed as $$\sup_{\epsilon>0}\limsup_{n\rightarrow \infty} \frac{1}{|F_n|}\log K_\varepsilon(f,X,G,\rho,F_n).$$

We use ideas in \cite[Section 5]{KerrLi10} to prove Theorem \ref{T-sofic pressure is equal to amenable pressure}. We need the following result, which is a Rokhlin lemma for sofic approximations \cite[Lemma 4.6]{KerrLi10}.
\begin{lemma}
\label{L-Rokhlin lemma for amenable group}
Let $G$ be a countable amenable group.
Let $0\le \tau<1$, $0<\eta<1$, $K$ be a nonempty finite subset of $G$, and $\delta>0$.
Then there are an $\ell\in \Nb$, nonempty finite sets $F_1, \dots, F_\ell \subset G$ with $\max_{1\leq k\leq \ell}\frac{|KF_k \setminus F_k|}{|F_k|}<\delta $ and
$\max_{1\leq k\leq \ell}\frac{|F_kK\setminus F_k|}{|F_k|}<\delta$, a finite subset $F$ of $G$ containing $e$,
and an $\eta'>0$ such that,
for every $d\in \Nb$, every map $\sigma: G\rightarrow \sym(d)$ for which there is a set $B\subset [d]$ satisfying
$|B|\geq (1-\eta')d$ and
$$
\sigma_s\sigma_t(a)=\sigma_{st}(a), \sigma_s(a)\neq \sigma_{s'}(a), \sigma_e(a)=a,
$$
for all $a\in B$ and $s, t, s'\in F$ with $s\neq s'$, and every set $V\subset [d]$ with $(1-\tau)d\leq|V|$,
there exist subsets $C_1, \dots, C_\ell$ of $V$ such that
\begin{enumerate}
\item for every $1\leq k\leq \ell$, the map $(s, c)\mapsto \sigma_s(c)$ from $F_k\times C_k$ to $\sigma(F_k)C_k$ is bijective,

\item the family $\{ \sigma(F_1)C_1, \dots, \sigma(F_\ell)C_\ell \}$ is disjoint and $(1-\tau-\eta)d\leq|\bigcup_{k=1}^\ell \sigma(F_k)C_k|$.
\end{enumerate}
\end{lemma}
\begin{lemma}
\label{L-sofic pressure is less or equal than amenable pressure}
Let $G$ be a countable amenable group acting continuously on a compact metrizable space $X$. Let $\Sigma$ be a sofic approximation sequence for $G$ and $f$ be a real valued continuous function on $X$. Then $P_{\Sigma}(f,X,G)\leq P(f,X,G)$.
\begin{proof}
We may assume that $P(f,X,G)<\infty$. Let $\rho$ be a compatible metric on $X$. It suffices to prove that $P_{\Sigma,\infty}(f,X,G,\rho)\leq P(f,X,G)+6\kappa$ for any $\kappa>0$.

Let $\kappa>0$. Let $\varepsilon'>0$ be such that $|f(x)-f(y)|<\kappa$ whenever $x,y\in X$ with $\rho(x,y)<\varepsilon'/2$. It suffices to prove that $P_{\Sigma,\infty}^\varepsilon(f,X,G,\rho)\leq P(f,X,G)+6\kappa$, for all $0<\varepsilon<\varepsilon'$.

Let $0<\varepsilon<\varepsilon'$. Then there are a nonempty finite set $K\subset G$ and $\delta'>0$ such that for any nonempty finite set
$F'\subset G$ satisfying $|KF'\setminus F'|<\delta'|F'|$, we have
$K_{\varepsilon/4}(f,X,G, \rho,F')<\exp((P(f,X,G) + \kappa)|F'|)$ . Since $f$ is continuous on $X$, there exists $Q>0$ such that $|f(x)|\leq Q$ for all $x\in X$.

Choose $0<\eta< 1$ such that $(N_{\varepsilon/4}(X,\rho))^{2\eta}\leq \exp(\kappa)$ and $\eta<\frac{\kappa}{2Q}$.

By Lemma \ref{L-Rokhlin lemma for amenable group} there are an $m\in \Nb$
and nonempty finite sets $F_{1}, \dots, F_{m}\subset G$
satisfying $\max_{1\leq k\leq m}\frac{|KF_k\setminus F_k|}{|F_k|}<\delta'$ such that for every good enough sofic approximation $\sigma : G\to\sym(d)$ for some $d\in \Nb$ and every $W\subset [d]$ with $(1-\eta)d\leq|W|$, there exist finite subsets $C_{1}, \dots, C_{m}$ of $W$ satisfying the following:
\begin{enumerate}
\item for every $k=1, \dots, m$, the map $(s, c)\mapsto \sigma_s(c)$ from $F_{k}\times C_{k}$ to $\sigma(F_{k})C_{k}$ is bijective,

\item the family $\{ \sigma(F_{1})C_{1}, \dots, \sigma(F_{m})C_{m} \}$ is disjoint and $(1-2\eta)d\leq|\bigcup_{k=1}^{m} \sigma(F_{k})C_{k}|$.
\end{enumerate}
Then $\max_{1\leq k\leq m}K_{\varepsilon/4}(f,X,G, \rho,F_k)\leq \exp\Big((P(f,X,G) + \kappa)|F_k|\Big).$

Let $\delta > 0$ and set $F=\bigcup_{k=1}^m F_k$.
Let $\sigma:G \rightarrow \sym(d)$ be a good enough sofic approximation of $G$, for some $d\in\Nb$.
We will show that $$M^{\varepsilon}_{\Sigma,\infty}(f,X,G,\rho,F, \delta, \sigma)\leq \exp((P(f,X,G) + 6\kappa)d),$$ when $\delta$ is small enough.

Let $\cE$ be a $(\rho_\infty,\varepsilon)$-separated subset of $\map(\rho,F,\delta,\sigma)$ such that $$M^{\varepsilon}_{\Sigma,\infty}(f,X,G,\rho,F,\delta,\sigma)\leq 2\cdot\sum_{\varphi\in \cE}\exp\Big(\sum_{a=1}^df(\varphi(a))\Big).$$

For each $\varphi\in \cE$ we denote by $\Lambda_\varphi$ the set of all $a\in[d]$ such that $\rho(\varphi(sa),s\varphi(a))<\sqrt{\delta}$ for all $s\in F$. Then $|\Lambda_\varphi|\geq (1-|F|\delta)d$. We enumerate the elements of $\{\Lambda_\varphi:\varphi\in \cE\}$ as $\Lambda_1,...,\Lambda_\ell$.
Then $\cE=\bigsqcup_{j=1}^\ell \cV_j$, where $\cV_j=\{\varphi\in \cE:\Lambda_\varphi=\Lambda_j\}$, for every $j=1,...,\ell$.

Choose $\delta>0$ such that $|F|\delta<\eta$ and $2\sqrt\delta<\varepsilon/4$. Then for any $j\in\{1,\dots, \ell\},$ there exist subsets $C_{j,1},...,C_{j,m}$ of $\Lambda_j$ such that
\begin{enumerate}
\item for every $1\leq k\leq m$, the map $(s, c)\mapsto \sigma_s(c)$ from $F_{k}\times C_{j,k}$ to $\sigma(F_{k})C_{j,k}$ is bijective,

\item the family $\{ \sigma(F_{1})C_{j,1}, \dots, \sigma(F_{m})C_{j,m} \}$ is disjoint and $(1-2\eta)d\leq |\bigcup_{k=1}^{m} \sigma(F_{k})C_{j,k}|$.
\end{enumerate}

Let $1\leq j\leq \ell,1\leq k\leq m$ and $c\in C_{j,k}$.
Let $\cW_{j,k,c}$ be a maximal $(\rho_{\sigma(F_k)c,\infty},\varepsilon/2)$-separated subset of $\cV_j$. Then $\cW_{j,k,c}$ is a $(\rho_{\sigma(F_k)c,\infty},\varepsilon/2)$-spanning subset of $\cV_j$.

For any two distinct elements $\varphi$ and $\psi$ of $\cW_{j,k,c}$ , since $c\in \Lambda_j= \Lambda_{\psi}=\Lambda_{\varphi}$, for every $s\in F_{k}$, we have
\begin{align*}
\rho(s\psi(c), s\varphi(c))&\ge \rho(\psi(sc), \varphi(sc))-\rho( \psi(sc), s\psi(c))-\rho(s\varphi(c), \varphi(sc))\\
&\ge  \rho(\psi(sc), \varphi(sc))-2\sqrt{\delta},
\end{align*}
and hence
\begin{align*}
\max_{s\in F_{k}}\rho(s\psi(c), s\varphi(c))\ge \max_{s\in F_{k}}\rho(\psi(sc), \varphi(sc))-2\sqrt{\delta}\geq\varepsilon/2-\varepsilon/4=\varepsilon/4.
\end{align*}
Thus $\{\varphi(c):\varphi\in \cW_{j,k,c}\}$ is a $(\rho_{F_k},\varepsilon/4)$-separated subset of $X$.

Choose $\delta>0$ such that $|f(x)-f(y)|<\kappa$ for all $x,y\in X$ with $\rho(x,y)<\sqrt\delta$. Then
\begin{eqnarray*}
&&\sum_{\varphi\in \cW_{j,k,c}}\exp\Big(\sum_{s\in F_k}f(\varphi(sc))\Big)\\
&=& \sum_{\varphi\in \cW_{j,k,c}}\exp\Big(\sum_{s\in F_k}f(s\varphi(c))\Big)\exp\Big(\sum_{s\in F_k}(f(\varphi(sc))-f(s\varphi(c)))\Big)\\
&\leq& \sum_{\varphi\in \cW_{j,k,c}}\exp\Big(\sum_{s\in F_k}f(s\varphi(c))\Big)\exp(|F_k|\kappa)\\
&\leq& K_{\varepsilon/4}(f,X,G,\rho,F_k)\exp(|F_k|\kappa)\\
&\leq& \exp\Big((P(f,X,G)+2\kappa)|F_k|\Big).
\end{eqnarray*}

  Let $\cW_j$ be a $(\rho_{\cZ_j,\infty},\varepsilon/2)$-spanning subset of $\cV_j$ with minimal cardinality, where $\cZ_j=[d]\setminus \bigcup_{k=1}^m\sigma(F_k)C_{j,k}$. Then
$$
|\cW_j | \leq (N_{\varepsilon/4}(X, \rho ))^{|\cZ_j|} \leq (N_{\varepsilon/4}(X, \rho ))^{2\eta d}\leq \exp(\kappa d).
$$

 Denote by $\cU_j$ the set of all maps $\varphi:[d]\rightarrow X$ such that $\varphi|_{\cZ_j}\in \cW_j|_{\cZ_j}$ and $\varphi|_{\sigma(F_k)c}\in \cW_{j,k,c}|_{\sigma(F_k)c}$ for all $1\leq k\leq m$ and $c\in C_{j,k}$. Then
\begin{eqnarray*}
&&\sum_{\varphi\in \cU_j}\exp\Big(\sum_{a=1}^df(\varphi(a))\Big)\\
&=&\sum_{\varphi\in \cU_j}\exp\Big(\sum_{k=1}^m\sum_{c\in C_{j,k}}\sum_{s\in F_k} f(\varphi(sc))\Big)\exp\Big(\sum_{a\in \cZ_j}f(\varphi(a))\Big)\\
&\leq&\sum_{\varphi\in \cU_j}\exp(2Q\eta d)\prod_{k=1}^m\prod_{c\in C_{j,k}}\exp\Big(\sum_{s\in F_k} f(\varphi(sc))\Big)\\
&\leq& (N_{\varepsilon/4}(X, \rho ))^{2\eta d}\exp(2Q\eta d)\prod_{k=1}^m\prod_{c\in C_{j,k}}\sum_{\psi\in \cW_{j,k,c}|_{\sigma(F_k)c}}\exp\Big(\sum_{s\in F_k} f(\psi(sc))\Big)\\
&\leq&(N_{\varepsilon/4}(X, \rho ))^{2\eta d}\exp(2Q\eta d)\prod_{k=1}^m\prod_{c\in C_{j,k}}\exp\Big((P(f,X,G)+2\kappa)|F_k|\Big)\\
&\leq& (N_{\varepsilon/4}(X, \rho ))^{2\eta d}\exp(2Q\eta d)\exp\Big((P(f,X,G)+2\kappa)\sum_{k=1}^m|F_k||C_{j,k}|\Big)
\\
&\leq& \exp(\kappa d)\exp(\kappa d)\exp\Big((P(f,X,G)+2\kappa)d\Big).
\end{eqnarray*}

 By spanning properties of $\cW_{j,k,c}$ and $\cW_j$, we can define a map $\Phi:\cV_j\rightarrow \cU_j$ by choosing for each $\psi\in \cV_j$, some $\Phi(\psi)\in \cU_j$ with $\rho_{\infty}(\psi,\Phi(\psi))\leq \varepsilon/2$. Then $\Phi$ is injective, so
 \begin{eqnarray*}
 \sum_{\psi\in \cU_j}\exp\Big(\sum_{a=1}^df(\psi(a))\Big)&\geq& \sum_{\psi\in \Phi(\cV_j)}\exp\Big(\sum_{a=1}^df(\psi(a))\Big)\\
 &=& \sum_{\varphi\in \cV_j}\exp\Big(\sum_{a=1}^d(f(\Phi(\varphi)(a))-f(\varphi(a)))\Big)\exp\Big(\sum_{a=1}^df(\varphi(a))\Big)\\
 &\geq& \exp(-d\kappa)\sum_{\varphi\in \cV_j}\exp\Big(\sum_{a=1}^df(\varphi(a))\Big).
 \end{eqnarray*}
Therefore
\begin{eqnarray*}
\sum_{\varphi\in \cE}\exp\Big(\sum_{a=1}^df(\varphi(a))\Big)&=&\sum_{j=1}^\ell\sum_{\varphi\in \cV_j}\exp\Big(\sum_{a=1}^df(\varphi(a))\Big)\\
&\leq&\sum_{j=1}^\ell\sum_{\varphi\in \cU_j}\exp\Big(\sum_{a=1}^df(\varphi(a))\Big)\exp(\kappa d)\\
&\leq& \ell \exp(\kappa d)\exp\Big((P(f,X,G)+2\kappa)d\Big)\exp(2\kappa d).
\end{eqnarray*}

By Stirling's approximation formula, $|F|\delta d \binom{d}{|F|\delta d}$ is less than $\exp(\beta d)$ for some $\beta>0$ depending on $\delta$ and $|F|$ but not on $d$ when $d$ is large enough with $\beta\rightarrow 0$ as $\delta\rightarrow 0$. Since $\sum_{j=0}^{\lfloor |F|\delta d\rfloor} \binom{d}{j}\leq |F|\delta d \binom{d}{|F|\delta d}$, when $d$ is large enough we have that the number of subsets of $[d]$ of cardinality at least $(1-|F|\delta)d$ is at most $\exp(\beta d)$. Choose $\delta$ such that $\beta<\kappa$. Then, when $d$ is large enough, $\ell\leq \exp(\beta d)\leq \exp(\kappa d)$.
Therefore
\begin{eqnarray*}
M^{\varepsilon}_{\Sigma,\infty}(f,X,G,\rho,F,\delta,\sigma)&\leq& 2\cdot\sum_{\varphi\in \cE}\exp\Big(\sum_{a=1}^df(\varphi(a))\Big)\\
&\leq& 2\cdot \exp(\kappa d)\exp(3\kappa d)\exp\Big((P(f,X,G)+2\kappa)d\Big),
\end{eqnarray*}
and hence $P^{\varepsilon}_{\Sigma,\infty}(f,X,G,\rho)\leq P^{\varepsilon}_{\Sigma,\infty}(f,X,G,\rho,F,\delta,\sigma)\leq P(f,X,G)+6\kappa,$ as we want.
\end{proof}
\end{lemma}
\begin{lemma}
\label{L-sofic pressure is greater than amenable pressure}
Let $G$ be a countable amenable group acting continuously on a compact metrizable space $X$ and $f$ a real valued continuous function on $X$. Then $P_{\Sigma}(f,X,G)\geq P(f,X,G)$.
\end{lemma}
\begin{proof}
Let $\rho$ be a compatible metric on $X$.

We will prove that for any real number $R< P(f,X,G)$ and $\kappa>0$, $P_{\Sigma,\infty}(f,X,G,\rho)\geq R-5\kappa$.
Let $R< P(f,X,G)$ and $\kappa>0$. Choose $\varepsilon_1>0$ such that $p_1(f,\varepsilon_1)>R-\kappa$.
Because $f$ is continuous, it is uniformly continuous on the compact space $X$. Thus, there exists $\varepsilon_2>0$ such that $|f(x)-f(y)|< \kappa$ for all $x,y\in X$ with $\rho(x,y)<\varepsilon_2$. Let $\varepsilon=\min\{\varepsilon_1,\varepsilon_2\}$.

For any nonempty finite subset $F'$ of $G$, and $(\rho_{F'},\varepsilon/2)$-separated subset $\cD$ of $X$ with maximal cardinality, $\{B_{F'}(x,\varepsilon/2)\}_{x\in \cD}$ is an open cover of $X$ of order $(F',\varepsilon)$, where $B_{F'}(x,\varepsilon/2)=\{y\in X: \max_{s\in F'}\rho(sx,sy)<\varepsilon/2\}$. Then $$|F'|^{-1}\log \sum_{x\in \cD}\sup_{y\in B_{F'}(x,\varepsilon/2)}\exp\Big (\sum_{s\in F'} f( sy)\Big)\geq p_1(f,\varepsilon)-\kappa,$$ whenever $F'$ is sufficiently left invariant.

 We also have $$\sum_{x\in \cD}\sup_{y\in B_{F'}(x,\varepsilon/2)}\exp\Big (\sum_{s\in F'} f(sy)\Big)\leq \exp (|F'|\kappa)\sum_{x\in \cD}\exp\Big  (\sum_{s\in F'} f(sx)\Big).$$
  Thus taking the logarithm on both sides, and dividing them by $|F'|$, when $F'$ is sufficiently left invariant, one has $$|F'|^{-1}\log\sum_{x\in \cD}\exp\Big (\sum_{s\in F'} f(sx)\Big)\geq p_1(f,\varepsilon)-2\kappa\geq R-3\kappa.$$
Then there exist a nonempty finite subset $K$ of $G$ and $\delta''>0$ such that $$\frac{1}{|F'|}\log\sum_{x\in \cD}\exp\Big (\sum_{s\in F'} f(sx)\Big)\geq R-3\kappa,$$
  for any nonempty finite subset $F'$ of $G$ satisfying $\frac{|KF'\setminus F'|}{|F'|}<\delta''$, and any $(\rho_{F'},\varepsilon/2)$-separated subset $\cD$ of $X$ with maximal cardinality.

Let $F$ be a nonempty finite subset of $G$ and $\delta>0$. We will show that if $\sigma:G\rightarrow \sym(d)$ is a good enough sofic approximation of $G$ then
\begin{align*} \label{E-top lower}
\frac{1}{d} \log M^{\varepsilon/2}_{\Sigma,\infty} (f,X,G,\rho ,F,\delta ,\sigma )\ge R - 5\kappa .
\end{align*}

Since $f$ is continuous on $X$ and $X$ is compact, there exists a number $Q>0$ such that $f(x)\geq -Q$ for all $x\in X$.
Choose $\delta'>0$ such that $\delta'<\delta'',\sqrt{\delta'}\diam(X,\rho)<\delta/\sqrt2$, $(1-\delta')(R-3\kappa)\geq R-4\kappa$ and $\delta'<\kappa/Q$.
By Lemma \ref{L-Rokhlin lemma for amenable group} there are an $\ell\in \Nb$
and nonempty finite sets $F_1, \dots, F_\ell\subset G$
satisfying $\max_{1\leq k\leq \ell}\frac{|KF_k\setminus F_k|}{|F_k|}<\delta'$ and $\max_{1\leq k\leq\ell}\frac{|FF_k\setminus F_k|}{|F_k|}<\delta'$
such that for every good
enough sofic approximation $\sigma : G\to\sym(d)$ for some $d\in \Nb$, and every $V\subset [d]$ with $(1-\delta'/2)d\leq|V|$, there exist subsets $C_1, \dots, C_\ell$ of $V$ satisfying the following:
\begin{enumerate}
\item for every $1\leq k\leq \ell$, the map $(t, c)\mapsto \sigma_t(c)$ from $F_k\times C_k$ to $\sigma(F_k)C_k$ is bijective,

\item the family $\{ \sigma(F_1)C_1, \dots, \sigma(F_\ell)C_\ell \}$ is disjoint and $(1-\delta')d\leq |\bigcup_{k=1}^\ell \sigma(F_k)C_k|$.
\end{enumerate}

For each map $\sigma: G\rightarrow \sym(d)$ for some $d\in\Nb$, put $$\Lambda_\sigma:=\{a\in [d]: \sigma_{st}(a)=\sigma_s\sigma_t(a) \mbox{ for any } s\in F \mbox{ and } t\in\bigcup _{k=1}^\ell F_k\}.$$
When $\sigma$ is a good enough approximation for $G$, one has $|\Lambda_\sigma|\geq (1-\delta'/2)d$. Then there exist $C_1,...,C_\ell\subset \Lambda_\sigma$ as above.

For each $1\leq k\leq \ell $, pick a $(\rho_{F_k},\varepsilon/2)$-separated subset
$\cE_k$ of $X$ with maximal cardinality. Then $$\frac{1}{|F_k|}\log\sum_{x\in \cE_k}\exp\Big (\sum_{s\in F_k} f(sx)\Big)\geq R-3\kappa,$$ for any $1\leq k\leq \ell$.

For every $h = (h_k )_{k=1}^\ell \in\prod_{k=1}^\ell (\cE_k )^{C_k}$
take a map $\varphi_h : [d]\rightarrow X$ such that
$$
\varphi_h(t c) = t(h_k (c))
$$
for all $1\leq k\leq \ell$, $t\in F_k$ and $c\in C_k$. Then for any $1\leq k\leq \ell$, $c\in C_k$, $s\in F$, and $t\in F_k$ satisfying $st\in F_k$, we have $\varphi_h(s(tc))=s\varphi_h(tc)$. Hence for any $s\in F$, one has
\begin{eqnarray*}
\sum_{k=1}^\ell\sum_{a\in \sigma(F_k)C_k}(\rho(\varphi_h(s(a)),s\varphi_h(a)))^2&=&\sum_{k=1}^\ell\sum_{c\in C_k}\sum_{t\in F_k, st\notin F_k}(\rho(\varphi_h(s(tc)),s\varphi_h(tc)))^2\\
&\leq & \sum_{k=1}^\ell|C_k||sF_k\setminus F_k|\diam^2(X,\rho)\\
&\leq & \sum_{k=1}^\ell|C_k||FF_k\setminus F_k|\diam^2(X,\rho)\\
&\leq & \sum_{k=1}^\ell|C_k||F_k|\delta'\diam^2(X,\rho)\\
&\leq & \delta'\diam^2(X,\rho)d.
\end{eqnarray*}
So
\begin{eqnarray*}
&&(\rho_2(\varphi_h\circ \sigma_s,\alpha_s\circ \varphi_h))^2\\
&=&\frac{1}{d}\Big(\sum_{k=1}^\ell\sum_{a\in \sigma(F_k)C_k}(\rho(\varphi_h(s(a)),s\varphi_h(a)))^2+\sum_{a\in [d]\setminus\bigcup_{k=1}^\ell\sigma(F_k)C_k}(\rho(\varphi_h(s(a)),s\varphi_h(a)))^2\Big)\\
&\leq& \delta'\diam^2(X,\rho)+\delta'\diam^2(X,\rho)< \delta,
\end{eqnarray*}
for any $s\in F$.
 Thus $\varphi_h \in\map (\rho , F ,\delta , \sigma )$.

For any distinct elements $h = (h_k )_{k=1}^\ell,h' = (h_k' )_{k=1}^\ell$ in
$\prod_{k=1}^\ell (\cE_k )^{C_k}$, there are a $1\leq k\leq \ell$ and a $c\in C_k$ such that $h_k(c)\neq h_k'(c)$.
Since $\cE_k$ is $(\rho_{F_k},\varepsilon/2)$-separated,
$\rho_{F_k}(h_k(c),h'_k(c))\geq \varepsilon/2$ and thus
we have
$\rho_{\infty} (\varphi_h ,\varphi_{h'} ) \ge \varepsilon/2$.
Then
\begin{eqnarray*}
M_{\Sigma,\infty}^{\varepsilon/2}(f,X,G,\rho,F,\delta,\sigma)&\geq& \sum_{h\in  \prod_{j=1}^\ell (\cE_j )^{C_j}}\exp\Big(\sum_{a=1}^df(\varphi_h(a))\Big)\\
&\geq& \sum_{h\in  \prod_{j=1}^\ell (\cE_j )^{C_j}}\exp \Big (\sum_{k=1}^l\sum_{c_k\in C_k}\sum_{s_k\in F_k}f(\varphi_h(s_kc_k))\Big )\exp(-Q\delta' d)\\
&=& \sum_{h\in  \prod_{j=1}^\ell (\cE_j )^{C_j}}\exp \Big (\sum_{k=1}^l\sum_{c_k\in C_k}\sum_{s_k\in F_k}f(s_kh(c_k))\Big )\exp(-Q\delta' d)\\
&=&\exp(-Q\delta' d)\sum_{h\in  \prod_{j=1}^\ell (\cE_j )^{C_j}}\prod_{k=1}^\ell\prod_{c_k\in C_k} \exp\Big(\sum_{s_k\in F_k}f(s_kh(c_k))\Big)\\
&=& \exp(-Q\delta' d)\prod_{j=1}^\ell \Big(\sum_{x\in \cE_j}\exp\big(\sum_{s\in F_j}f(sx)\big)\Big)^{|C_j|}.
\end{eqnarray*}
Therefore,
\begin{eqnarray*}
\frac{1}{d}\log M_{\Sigma,\infty}^{\varepsilon/2}(f,X,G,\rho,F,\delta,\sigma)&\geq&\frac{1}{d}\log \prod_{j=1}^\ell \Big(\sum_{x\in \cE_j}\exp\big(\sum_{s\in F_j}f(sx)\big)\Big)^{|C_j|}-Q\delta'\\
&=&\frac{1}{d}\sum_{j=1}^\ell|C_j|\log \Big(\sum_{x\in \cE_j}\exp\big(\sum_{s\in F_j}f(sx)\big)\Big)-Q\delta'\\
&\geq& \frac{1}{d}\sum_{j=1}^\ell(R-3\kappa)|C_j||F_j|-\kappa.
\end{eqnarray*}

If $R-3\kappa\geq 0$ then $\frac{1}{d}\sum_{j=1}^\ell(R-3\kappa)|C_j||F_j|\geq (1-\delta')(R-3\kappa)\geq R-4\kappa$ and if $R-3\kappa< 0$ then $\frac{1}{d}\sum_{j=1}^\ell(R-3\kappa)|C_j||F_j|\geq R-3\kappa\geq R-4\kappa$. Thus, $\frac{1}{d}\log M_{\Sigma,\infty}^{\varepsilon/2}(f,X,G,\rho,F,\delta,\sigma)\geq R-5\kappa$, as desired.
\end{proof}
Combining Lemmas \ref{L-sofic pressure is less or equal than amenable pressure} and \ref{L-sofic pressure is greater than amenable pressure} we obtain Theorem 1.1.
\section{The variational principle of topological pressure}
We will prove Theorem 1.2 in this section. Let $\alpha$ be a continuous action of a countable sofic group $G$ on a compact metrizable space $X$. Before proving the variational principle for sofic topological pressure, we recall the definition of sofic measure entropy \cite[Section 3]{KerrLi10}.
\subsection{Sofic measure entropy}
Let $\mu$ be a Borel probability measure on $X$ and $\rho$ a continuous pseudometric on $X$.
\begin{definition}
Let $L$ be a nonempty finite subset of $C(X)$, $F$ a nonempty finite subset of $G$,  and $\delta > 0$.
Let $\sigma$ be a map from $G$ to $\sym (d)$ for some $d\in\Nb$.
We define $\map_\mu (\rho ,F, L,\delta ,\sigma )$ to be the set of all $\varphi$ in $\map(\rho ,F,\delta ,\sigma )$ such that
$$\big|\frac{1}{d}\sum_{j=1}^df(\varphi(j))-\int_X f\,d\mu\big| < \delta, \mbox { for all } f\in L.$$
\end{definition}
\begin{definition}
For $\varepsilon > 0$ we define
\begin{align*}
h_{\Sigma ,\mu ,\infty}^\varepsilon (\rho ,F, L,\delta ) &=
\limsup_{i\to\infty} \frac{1}{d_i} \log N_\varepsilon (\map_\mu (\rho ,F, L,\delta ,\sigma_i ),\rho_\infty ) ,\\
h_{\Sigma ,\mu ,\infty}^\varepsilon (\rho ,F, L) &= \inf_{\delta > 0} h_{\Sigma ,\mu ,\infty}^\varepsilon (\rho ,F,L,\delta ) ,\\
h_{\Sigma ,\mu ,\infty}^\varepsilon (\rho ,F) &= \inf_{L} h_{\Sigma ,\mu ,\infty}^\varepsilon (\rho ,F,L) ,\\
h_{\Sigma ,\mu ,\infty}^\varepsilon (\rho ) &= \inf_{F} h_{\Sigma ,\mu ,\infty}^\varepsilon (\rho ,F) ,\\
h_{\Sigma ,\mu ,\infty} (\rho ) &= \sup_{\varepsilon > 0} h_{\Sigma ,\mu ,\infty}^\varepsilon (\rho ) ,
\end{align*}
where $L$ in the third line runs over the nonempty finite subsets of $C(X)$ and
$F$ in the fourth line runs over the nonempty finite subsets of $G$.

If $\map_\mu (\rho ,F, L,\delta ,\sigma_i )=\emptyset$ for all
large enough $i$, we set $h_{\Sigma ,\mu ,\infty}^\varepsilon (\rho ,F, L,\delta ) = -\infty$.
\end{definition}
If $\mu$ is a $G$-invariant Borel probability measure on $X$ and $\rho$ is a dynamically generating pseudometric then from Proposition 5.4 in \cite{KerrLi11a} and Proposition 3.4 in \cite{KerrLi10}, we conclude that $h_{\Sigma ,\mu ,\infty} (\rho )$ coincides with the sofic measure entropy $h_{\Sigma,\mu}(X,G)$ (see \cite{KerrLi11a} for the definition of $h_{\Sigma,\mu}(X,G)$). In particular, the quantities $h_{\Sigma ,\mu ,\infty} (\rho )$ do not depend on the choice of compatible metrics on $X$.

Now we prove the variational principle for sofic topological pressure.
\subsection{The variational principle}
We denote by $M(X)$ the convex set of Borel probability measures on $X$. Denote by $M_G(X)$ the set of $G$-invariant Borel probability measures on $X$. Under the weak* topology, $M(X)$ is compact and $M_G(X)$ is a closed convex subset of $M(X)$.

The following Lemma was proved by Kerr and Li in \cite[Theorem 6.1]{KerrLi11a} for the case $f=0$. We modify the argument there to deal with general functions $f$ in $C(X)$.
\begin{lemma}
\label{L-sofic pressure is less than or equal sum of measure entropy and mu(f)}
Let $\alpha$ be a continuous action of a countable sofic group $G$ on a compact metrizable space $X$. Let $\Sigma$ be a sofic approximation sequence for $G$ and $f$ be a real valued continuous function on $X$. Then $$P_{\Sigma,\infty}(f,X,G)\leq\sup\Big\{h_{\Sigma,\mu}(X,G)+\int_X f\,d\mu:\mu\in M_G(X)\Big\}.$$
\end{lemma}
\begin{proof}
Let $\rho$ be a compatible metric on $X$. We may assume that $P_{\Sigma,\infty}(f,X,G)\neq -\infty$. Let $\varepsilon>0$. It suffices to prove that there exists $\mu\in M_G(X)$ such that $h^\varepsilon_{\Sigma,\mu,\infty}(\rho)+\int_X fd\mu\geq P^\varepsilon_{\Sigma,\infty}(f,X,G,\rho)$.

Take a sequence $e\in F_1\subset F_2\subset\dots$ of finite subsets of $G$ such that $G=\bigcup_{n\in\Nb}F_n$. Since $X$ is compact and metrizable, there exists a sequence $\{g_m\}_{m\in \Nb}$ in $C(X)$ such that $\{g_m\}_{m\in \Nb}$ is dense in $C(X)$. Let $n\in \Nb$ and $L_n=\{f,g_1,\dots g_n\}$. There exists $Q>0$ such that $\max_{g\in L_n}\|g\|_\infty\leq Q$.
Choose $\delta_n>0$ such that $\delta_n< \frac{1}{12Q|F_n|},\delta_n<\frac{1}{3n}$ and $|g(x)-g(y)|<\frac{1}{6n}$ for all $g\in L_n$ and for all $x,y\in X$ with $\rho(x,y)<\sqrt{\delta_n}$.
We will find some $\mu_n\in M(X)$ such that $$h^\varepsilon_{\Sigma,\mu_n,\infty}(\rho,F_n,L_n,\frac{1}{3n})+\int_X fd\mu_n+\frac{1}{3n}\geq P^\varepsilon_{\Sigma,\infty}(f,X,G,\rho),$$ and $|\mu_n(\alpha_{t^{-1}}(g))-\mu_n(g)|<1/n$ for any $t\in F_n$, $g\in L_n$.

Since $M(X)$ is compact under weak* topology, there exists a finite subset $\cD$ of $M(X)$ such that for any map $\sigma:G\rightarrow \sym(d)$ for some $d\in \Nb$ and any $\varphi\in \map(\rho,F_n,\delta_n,\sigma)$ there is a $\mu_\varphi\in\cD$ such that $|\mu_\varphi(\alpha_{t^{-1}}(g))-(\varphi_*\zeta)(\alpha_{t^{-1}}(g))|<\frac{1}{3n}$ for all $t\in F_n$, $g\in L_n$, where $\zeta$ is the uniform probability measure on $[d]$, i.e., $(\varphi_*\zeta)(h)=\frac{1}{d}\sum_{a=1}^d h(\varphi(a))$ for all $h\in C(X)$.

Let $\sigma$ be a map from $G$ to $\sym(d)$ for some $d\in\Nb$. For each $\varphi\in\map(\rho,F_n,\delta_n,\sigma)$, denote by $\Lambda_\varphi$ the set of all $a$ in $[d]$ such that $\rho(\varphi(ta),t\varphi(a))<\sqrt{\delta_n}$ for all $t\in F_n$. Then $|\Lambda_\varphi|\geq (1-|F_n|\delta_n)d.$ Thus, for all $t\in F_n$, $g\in L_n$, we have
\begin{eqnarray*}
|(\varphi_*\zeta)(\alpha_{t^{-1}}(g))-((\varphi\circ \sigma_{t})_*\zeta)(g)|&\leq&\frac{1}{d}\Big|\sum_{a\in\Lambda_\varphi}(g(t\varphi(a))-g(\varphi(ta)))\Big|\\
&&+\frac{1}{d}\Big|\sum_{a\notin\Lambda_\varphi}(g(t\varphi(a))-g(\varphi(ta)))\Big|\\
&\leq&\frac{1}{d}|\Lambda_\varphi|\cdot \frac{1}{6n}+\frac{1}{d}2Q|F_n|\delta_nd \\
&\leq& \frac{1}{6n}+\frac{1}{6n}=\frac{1}{3n},
\end{eqnarray*}
and hence
\begin{eqnarray*}
|\mu_\varphi(\alpha_{t^{-1}}(g))-\mu_\varphi(g)|&\leq& |\mu_\varphi(\alpha_{t^{-1}}(g))-(\varphi_*\zeta)(\alpha_{t^{-1}}(g))|+|(\varphi_*\zeta)(g)-\mu_\varphi(g)|\\
&&+|(\varphi_*\zeta)(\alpha_{t^{-1}}(g))-((\varphi\circ \sigma_{t})_*\zeta)(g)|\\
&\leq& \frac{1}{3n}+\frac{1}{3n}+\frac{1}{3n}=\frac{1}{n}.
\end{eqnarray*}

Take a maximal $(\rho_\infty,\varepsilon)$-separated subset $\cE_\sigma$ of $\map(\rho, F_n,\delta_n,\sigma)$ such that $$M^{\varepsilon}_{\Sigma,\infty}(f,X,G,\rho,F_n,\delta_n,\sigma)\leq \exp(1)\cdot\sum_{\varphi\in \cE_\sigma}\exp\Big(\sum_{a=1}^df(\varphi(a))\Big).$$
For any $\nu\in\cD$, we denote by $W(\sigma,\nu)$ the set of all elements $\varphi$ in $\cE_\sigma$ such that $\mu_\varphi=\nu$.
By the pigeonhole principle there exists a $\nu_0\in\cD$ such that
$$|\cD|\cdot\sum_{\varphi\in W(\sigma,\nu_0)}\exp\Big({\sum_{a=1}^df(\varphi(a))}\Big)\geq \sum_{\varphi\in \cE_\sigma}\exp\Big(\sum_{a=1}^df(\varphi(a))\Big).$$
Since $|\nu_0(f)-(\varphi_*\zeta)(f)|<\frac{1}{3n}$ for all $\varphi\in W(\sigma,\nu_0)$, we have $\exp(\nu_0(f)d+\frac{d}{3n})\geq \exp\Big(\sum_{a=1}^df(\varphi(a))\Big)$ for all $\varphi\in W(\sigma,\nu_0)$ and hence
\begin{eqnarray*}
|\cD||\cW(\sigma,\nu_0)|\exp(\nu_0(f)d+\frac{d}{3n})&\geq& |\cD|\cdot\sum_{\varphi\in \cW(\sigma,\nu_0)}\exp\Big(\sum_{a=1}^df(\varphi(a))\Big)\\
&\geq&\sum_{\varphi\in \cE_\sigma}\exp\Big(\sum_{a=1}^df(\varphi(a))\Big).
\end{eqnarray*}
Note that $\cW(\sigma,\nu_0)\subset \map_{\nu_0}(\rho,F_n,L_n,\frac{1}{3n},\sigma)$ as $e\in F_n$ and $\delta_n<\frac{1}{3n}$. Since $\cW(\sigma,\nu_0)$ is $(\rho_\infty,\varepsilon)$-separated, we obtain
\begin{eqnarray*}
\frac{1}{d}\log \sum_{\varphi\in \cE_\sigma}\exp\Big(\sum_{a=1}^df(\varphi(a))\Big)&\leq&\frac{1}{d}\log(|\cD||\cW(\sigma,\nu_0)|)+\nu_0(f)+\frac{1}{3n}\\
&\leq&\frac{1}{d}\log\Big(|\cD|N_\varepsilon(\map_{\nu_0}(\rho,F_n,L_n,\frac{1}{3n},\sigma))\Big)+\nu_0(f)+\frac{1}{3n}.
\end{eqnarray*}
Thus
\begin{eqnarray*}
&&\frac{1}{d}\log M^\varepsilon_{\Sigma,\infty}(f,X,G,\rho, F_n,\delta_n,\sigma)\\
&\leq&\frac{1}{d}+\frac{1}{d}\log\Big(\sum_{\varphi\in \cE_{\sigma}}\exp\Big(\sum_{a=1}^{d}f(\varphi(a))\Big)\Big)\\
&\leq& \frac{1}{d}+\frac{1}{d}\log\Big(|\cD|N_\varepsilon(\map_{\nu_0}(\rho,F_n,L_n,\frac{1}{3n},\sigma))\Big)
+\nu_0(f)+\frac{1}{3n}.
\end{eqnarray*}

Letting $\sigma$ run through the terms of the sofic approximation sequence $\Sigma$, by the pigeonhole principle there exist $\mu_n\in \cD$ and a sequence $i_1<i_2<\dots $ in $\Nb$ with $$P^\varepsilon_{\Sigma,\infty}(f,X,G,\rho,F_n,\delta_n)=\lim_{k\rightarrow\infty}\frac{1}{d_{i_k}}\log M^\varepsilon_{\Sigma,\infty}(f,X,G,\rho, F_n,\delta_n,\sigma_{i_k})$$ such that
\begin{eqnarray*}
\frac{1}{d_{i_k}}\log M^\varepsilon_{\Sigma,\infty}(f,X,G,\rho, F_n,\delta_n,\sigma)&\leq& \frac{1}{d_{i_k}}\log\Big(|\cD|N_\varepsilon(\map_{\mu_n}(\rho,F_n,L_n,\frac{1}{3n},\sigma_{i_k}))\Big)\\
&&+\frac{1}{d_{i_k}}+\mu_n(f)+\frac{1}{3n},
\end{eqnarray*}
for all $k\in\Nb$ and $|\mu_n(\alpha_{t^{-1}}(g))-\mu_n(g)|<1/n$ for any $t\in F_n$, and $g\in L_n$. Then
\begin{eqnarray*}
&&P^\varepsilon_{\Sigma,\infty}(f,X,G,\rho)\\
&\leq& P^\varepsilon_{\Sigma,\infty}(f,X,G,\rho,F_n,\delta_n)\\
&=& \lim_{k\rightarrow\infty}\frac{1}{d_{i_k}}\log M^\varepsilon_{\Sigma,\infty}(f,X,G,\rho, F_n,\delta_n,\sigma_{i_k})\\
&\leq& \lim_{k\rightarrow\infty}\Big(\frac{1}{d_{i_k}}+\frac{1}{d_{i_k}}\log\Big(|\cD|N_\varepsilon(\map_{\mu_n}(\rho,F_n,L_n,\frac{1}{3n},\sigma_{i_k}))\Big)+\mu_n(f)+\frac{1}{3n}\Big)\\
&\leq& h^\varepsilon_{\Sigma,\mu_n,\infty}(\rho,F_n, L_n,\frac{1}{3n})+\mu_n(f)+\frac{1}{3n}.
\end{eqnarray*}

Let $\mu$ be a weak* limit point of the sequence $\{\mu_n\}_{n=1}^\infty$. For any $t\in G$ and $g\in\{g_m\}_{m\in\Nb}$, we have
$$|\mu(\alpha_{t^{-1}}(g))-\mu(g)|\leq |\mu(\alpha_{t^{-1}}(g))-\mu_n(\alpha_{t^{-1}}(g))|+|\mu_n(\alpha_{t^{-1}}(g))-\mu_n(g)|+|\mu_n(g)-\mu(g)|.$$
Since the right hand side converges to 0 as $n\rightarrow \infty$ and $\{g_m\}_{m\in\Nb}$ is dense in $C(X)$, we deduce that $\mu$ is $G$-invariant.

Let $F$ be a nonempty finite subset of $G$, $L$ a nonempty finite subset of $C(X)$ and $\delta>0$. Choose $n\in \Nb$ such that $F\subset F_n, \frac{1}{3n}\leq \delta/4$, $\max_{g\in L\cup \{f\}}|\mu_n(g)-\mu(g)|<\delta/4$ and for any $g\in L$, there exists $g'\in L_n$ such that $\|g-g'\|_\infty<\delta/4$. 
Then for any map $\sigma:G \rightarrow \sym(d)$ for some $d\in \Nb$, $\varphi\in \map_{\mu_n}(\rho,F_n,L_n,\frac{1}{3n},\sigma)$ and $g\in L$, we have
\begin{eqnarray*}
|(\varphi_*\zeta)(g)-\mu(g)|&\leq & |(\varphi_*\zeta)(g)-(\varphi_*\zeta)(g')|+|(\varphi_*\zeta)(g')-\mu_n(g')|\\
&&+|\mu_n(g')-\mu_n(g)|+|\mu_n(g)-\mu(g)|\\
&<& \frac{3\delta}{4}+\frac{1}{3n}\leq \delta,
\end{eqnarray*}
and hence $\varphi\in \map_\mu(\rho, F,L,\delta,\sigma)$. Thus $$\map_{\mu_n}(\rho, F_n,L_n,\frac{1}{3n},\sigma)\subset \map_\mu(\rho, F,L,\delta,\sigma)$$
and then
\begin{eqnarray*}
h^\varepsilon_{\Sigma,\mu,\infty}(\rho,F,L,\delta)+\int_X fd\mu
&\geq& h^\varepsilon_{\Sigma,\mu_n,\infty}(\rho,F_n,L_n,\frac{1}{3n})+\int_X fd\mu_n-\frac{\delta}{4}\\
&\geq& P^\varepsilon_{\Sigma,\infty}(f,X,G,\rho)-\frac{1}{3n}-\frac{\delta}{4}\\
&\geq& P^\varepsilon_{\Sigma,\infty}(f,X,G,\rho)-\frac{\delta}{2} .
\end{eqnarray*}
 Since $F,L,\delta$ are arbitrary we get $h^\varepsilon_{\Sigma,\mu,\infty}(\rho)+\int_X fd\mu\geq P^\varepsilon_{\Sigma,\infty}(f,X,G,\rho),$ as desired. Then $$P_{\Sigma,\infty}(f,X,G)\leq\sup\Big\{h_{\Sigma,\mu}(X,G)+\int_X f\,d\mu:\mu\in M_G(X)\Big\}.$$
\end{proof}
We can now prove Theorem \ref{T-variational principle for topological pressure}.
\begin{proof}[Proof of Theorem 1.2.]
Let $\rho$ be a compatible metric on $X$ and $\mu\in M_G(X)$. Let $F$ be a nonempty finite subset of $G$, and $\delta,\varepsilon>0$. Put $L_1=\{f\} $. Fix $i\in\Nb$. Let $\cE$ be a $(\rho_\infty,\varepsilon)$-separated subset of $\map_\mu(\rho,F,L_1,\delta,\sigma_i)$ with maximal cardinality. Then $\cE$ is also a $(\rho_\infty,\varepsilon)$-separated subset of $\map(\rho,F,\delta,\sigma_i)$.

Since the function $x\mapsto \log x$ for $x>0$ is concave, one has $$\log\sum_{\varphi\in \cE}\frac{1}{|\cE|}\exp\Big(\sum_{j=1}^{d_i}f(\varphi(j))\Big)\geq \frac{1}{|\cE|} \sum_{\varphi\in \cE}\sum_{j=1}^{d_i}f(\varphi(j)).$$
Hence
 \begin{eqnarray*}
\log\sum_{\varphi\in \cE}\exp\Big(\sum_{j=1}^{d_i}f(\varphi(j))\Big)&\geq& \log |\cE| +\frac{1}{|\cE|} \sum_{\varphi\in \cE}\sum_{j=1}^{d_i}f(\varphi(j))\\
&\geq& \log |\cE|+\frac{1}{|\cE|}\sum_{\varphi\in \cE}\Big(\int_X f\,d\mu-\delta\Big)d_i\\
&=& \log |\cE|+\Big(\int_X f\,d\mu-\delta\Big)d_i.
\end{eqnarray*}
Thus $P^\varepsilon_{\Sigma,\infty}(f,X,G,\rho,F,\delta)+\delta\geq h^\varepsilon_{\Sigma,\mu,\infty}(\rho,F,L_1,\delta)+\int_X f\,d\mu,$ for all nonempty finite subset $F$ of $G$ and all $\delta,\varepsilon>0$, yielding $P^\varepsilon_{\Sigma,\infty}(f,X,G,\rho,F)\geq h^\varepsilon_{\Sigma,\mu,\infty}(\rho,F,L_1)+\int_X f\,d\mu\geq h^\varepsilon_{\Sigma,\mu,\infty}(\rho,F)+\int_X f\,d\mu$ for all nonempty finite subset $F$ of $G$ and any $\varepsilon>0$. Hence $P_{\Sigma}(f,X,G)\geq h_{\Sigma,\mu}(X,G)+\int_X f\,d\mu.$
Combining with Lemma \ref{L-sofic pressure is less than or equal sum of measure entropy and mu(f)}, we get
$$P_{\Sigma}(f,X,G)=\sup\Big\{h_{\Sigma,\mu}(X,G)+\int_X fd\mu:\mu\in M_G(X)\Big\}.$$
\end{proof}
\begin{remark}
From the variational principle theorem we see that if $X$ has no $G$-invariant Borel probability measure then the topological pressure will be $-\infty$. For an example of such action, see the example at the end of section 4 in \cite{KerrLi11a}.
Note that when $G$ is amenable, for any continuous action of $G$ on a compact metrizable space there always exists a $G$-invariant Borel probability measure. In this case, the sofic topological pressure is always different from $-\infty$ since it coincides with the classical topological pressure; see Theorem \ref{T-sofic pressure is equal to amenable pressure}.
\end{remark}

\section{Equilibrium States and Examples}
 In this section we will calculate the sofic topological pressure of some functions over a Bernoulli shift. Let $\alpha$ be a continuous action of a countable sofic group $G$ on a compact metrizable space $X$.
\begin{definition}
Let $\Sigma$ be a sofic approximation sequence of $G$ and $f$ be a real valued continuous function on $X$.
A member $\mu$ of $M_G(X)$ is called an \textit{equilibrium state for f} with respect to $\Sigma$ if $P_\Sigma(f,X,G)=h_{\Sigma ,\mu}(X,G)+\int_X fd\mu$.
\end{definition}
\begin{definition}
Let $Y=\{0,...,k-1\}$ for some $k\in\Nb$ and $\mu$ a probability measure on $Y$. Let $Y^G=\prod_{s\in G}Y$ be the set of all functions $y:G\rightarrow Y$. For any nonempty finite subset $F$ of $G$, $a=(a_s)_{s\in F}\in Y^F$, put $A_{F,a}=\{(y_t)_{t\in G}: y_s=a_s, \mbox{ for all }s\in F\}$. Then there exists a unique measure $\mu^G$ on $Y^G$ defined on the $\sigma$-algebra of Borel subsets of $Y^G$ such that $\mu^G(A_{F,a})=\prod_{s\in F}\mu(a_s)$ for any nonempty finite subset $F$ of $G$, and $a=(a_s)_{s\in F}\in Y^F$, see \cite[page 5]{Walters82}.
\end{definition}

The following result is known when the acting group $G=\Zb^d$ for some $d\in \Nb$. For example, see \cite[Theorem 9.16]{Walters82} for the case $d=1$ and \cite[Example 4.2.2]{Keller98} for the general case $d\in \Nb$.
\begin{theorem}
\label{T-example}
Let $G$ be a countable sofic group, $k\in\Nb$ and $X=\{0,1,...,k-1\}^G$. Let $a_0,...,a_{k-1}\in\Rb$ and define $f\in C(X)$ by $f(x)=a_{x_e}$ where $x=(x_t)_{t\in G}$. Let $\alpha$ be the continuous action of $G$ on $X^G$ by the left shifts $s\cdot (x_t)_{t\in G}=(x_{s^{-1}t})_{t\in G}$. Let $\Sigma$ be a sofic approximation sequence of $G$ and $\mu$ the probability measure on $\{0,...,k-1\}$, defined by $$\mu(i)=\frac{\exp(a_i)}{\sum_{j=0}^{k-1}\exp(a_j)}, \mbox{   for all }0\leq i\leq k-1 .$$
Then
\begin{eqnarray*}
P_{\Sigma}(f,X,G)&=&\sup
\Big\{H(p)+\sum_{i=0}^{k-1}p(i)a_i: \mbox{ p is a probability measure on }\{0,...,k-1\}\Big\}\\
&=&\log \Big(\sum_{j=0}^{k-1}\exp(a_j)\Big),
\end{eqnarray*}
 where $H(p)=\sum_{i=0}^{k-1}-p(i)\log p(i)$.
Furthermore, the measure $\mu^G$ is an equilibrium state for $f$.
\end{theorem}
\begin{proof}
 Let $\rho$ be the pseudometric on $X$ defined by $\rho(x,y)=1$ if $x_e\neq y_e$ and $\rho(x,y)=0$ if $x_e=y_e$, where $x=(x_s)_{s\in G}, y=(y_s)_{s\in G}\in X$. Then $\rho$ is a continuous dynamically generating pseudometric on $X$. Let $1>\varepsilon>0,\delta>0$ and $F$ be a nonempty finite subset of $G$. Let $\sigma$ be a map from $G$ to $\sym(d)$ for some $d\in\Nb$. Let $\cE$ be a $(\rho_\infty,\varepsilon)$-separated subset of $\map(\rho,F,\delta,\sigma)$. Since $\cE$ is $(\rho_\infty,\varepsilon)$-separated, for any distinct elements $\varphi,\psi\in\cE$, $(\varphi(j))_e\neq (\psi(j))_e$ for some $1\leq j\leq d$. Thus
\begin{eqnarray*}
\sum_{\varphi\in\cE}\exp\Big(\sum_{j=1}^df(\varphi(j))\Big)&=&\sum_{\varphi\in\cE}\exp\Big(\sum_{j=1}^da_{(\varphi(j))_e}\Big)\\
&\leq&\sum_{(b_1,...,b_d)\in\{a_0,...,a_{k-1}\}^d}\exp\Big(\sum_{j=1}^db_j\Big)\\
&=&\sum_{(b_1,...,b_d)\in\{a_0,...,a_{k-1}\}^d}\prod_{j=1}^d\exp(b_j)\\
&=&\Big(\sum_{i=0}^{k-1}\exp(a_i)\Big)^d,
\end{eqnarray*}
and hence $\frac{1}{d}\log M^\varepsilon_{\Sigma,\infty}(f,X,G,F,\delta,\sigma)\leq \log \Big(\sum_{i=0}^{k-1}\exp(a_i)\Big)$.

For each $\beta\in\{0,...,k-1\}^d$, take a map $\varphi_\beta:\{1,...,d\}\rightarrow X^G$ such that for each $i\in [d]$ and $t\in G$, $((\varphi_\beta)(i))_t=\beta(\sigma(t^{-1})i)$. We denote by $\cZ$ the set of $i$ in $[d]$ such that $\sigma(e)\sigma(s)i=\sigma(s)i$ for all $s\in F$. For every $\beta\in \{0,...,k-1\}^d,s\in F$ and $i\in\cZ$, we have $(s\varphi_\beta(i))_e=(\varphi_\beta(i))_{s^{-1}}=\beta(\sigma(s)i)$ and $(\varphi_\beta(si))_e=\beta(\sigma(e)si)$, and hence $(s\varphi_\beta(i))_e=(\varphi_\beta(si))_e$.

When $\sigma$ is a good enough sofic approximation of $G$, one has $1-|\cZ|/d<\delta^2$, and hence $\varphi_\beta\in \map(\rho,F,\delta,\sigma)$. Note that $\{\varphi_\beta\}_{\beta\in \{0,...,k-1\}^d}$ is $(\rho_\infty,\varepsilon)$-separated. Thus
\begin{eqnarray*}
\frac{1}{d}\log M^\varepsilon_{\Sigma,\infty}(f,X,G,F,\delta,\sigma)&\geq& \frac{1}{d}\log \sum_{\beta\in \{0,...,k-1\}^d}\exp\Big(\sum_{i=1}^df(\varphi_\beta(i))\Big)\\
&=&\frac{1}{d}\log \sum_{\beta\in \{0,...,k-1\}^d}\exp\Big(\sum_{i=1}^da_{(\varphi_\beta(i))_e}\Big)\\
&=&\frac{1}{d}\log \sum_{\beta\in \{0,...,k-1\}^d}\exp\Big(\sum_{i=1}^da_{\beta(\sigma(e)i)}\Big)\\
&=&\frac{1}{d}\log \Big(\sum_{i=0}^{k-1}\exp(a_i)\Big)^d\\
&=&\log \Big(\sum_{i=0}^{k-1}\exp(a_i)\Big),
\end{eqnarray*}
and hence $\frac{1}{d}\log M^\varepsilon_{\Sigma,\infty}(f,X,G,F,\delta,\sigma)=\log \Big(\sum_{i=0}^{k-1}\exp(a_i)\Big)$. Thus $$P_\Sigma(f,X,G)=\log \Big(\sum_{i=0}^{k-1}\exp(a_i)\Big).$$

Let $\nu\in M_G(X)$. Put $A_i=\{(x_s)_{s\in G}\in X: x_e=i\}$ for any $i=0,...,k-1$. Let $p$ be the probability measure on $\{0,...,k-1\}$, defined by $p(i)=\nu(A_i)$ for any $i=0,...,k-1$. Then $$\int_X fd\nu=\sum_{i=0}^{k-1}\int_{A_i}f d\nu=\sum_{i=0}^{k-1}a_i\nu(A_i)=\sum_{i=0}^{k-1}a_ip(i)=\int_X fdp^G.$$
Since $\xi=\{A_0,...,A_{k-1}\}$ is a finite generating measurable partition of $X$, applying \cite[Proposition 5.3]{Bowen10} (taking $\beta$ there to be the trivial partition), \cite[Theorem 3.6]{KerrLi11a} and \cite[Proposition 3.4]{KerrLi10}, we get $h_{\Sigma,\nu}(X,G)\leq H_\nu(\xi)$, where $H_\nu(\xi)=\sum_{i=0}^{k-1}-\nu(A_i)\log \nu(A_i)$. Hence by Lemma 9.9 of \cite{Walters82},
\begin{eqnarray*}
h_{\Sigma,\nu}(X,G)+\int_X fd\nu
&\leq&H_\nu(\xi)+\sum_{i=0}^{k-1}a_ip(i)\\
&=&\sum_{i=0}^{k-1}p(i)(a_i-\log p(i))\\
&\leq&\log \Big(\sum_{i=0}^{k-1}\exp(a_i)\Big),
\end{eqnarray*}
From combining \cite[Theorem 8.1]{Bowen10}, \cite[Theorem 3.6]{KerrLi11a} and \cite[Proposition 3.4]{KerrLi10}, we know that the inequality in the first line becomes equality when $\nu=p^G$.
Furthermore, by Lemma 9.9 of \cite{Walters82}, the inequality in the third line becomes equality iff $$p(i)=\frac{\exp(a_i)}{\sum_{j=0}^{k-1}\exp(a_j)}=\mu(i), \mbox{ for every } 0\leq i\leq k-1.$$
Thus
\begin{eqnarray*}
P_{\Sigma}(f,X,G)&=&\sup
\Big\{H(p)+\sum_{i=0}^{k-1}p(i)a_i: \mbox{ p is a probability measure on }\{0,...,k-1\}\Big\}\\
&=&\log \Big(\sum_{j=0}^{k-1}\exp(a_j)\Big),
\end{eqnarray*}
and $\mu^G$ is an equilibrium state for $f$.
\end{proof}
When $G=\Zb$, $\mu^G$ is the unique equilibrium state for $f$, for example, see \cite[Theorem 9.16]{Walters82}. The proof there also works for the case $G$ is countable amenable. Thus, we raise the following question.
\begin{question}
Let $G$ be a countable sofic group, $k\in \Nb$ and $X, f\in C(X), \alpha,\mu$ as in the assumptions of Theorem \ref{T-example}. Is $\mu^G$ the unique equilibrium state for $f$ with respect to $\Sigma$, for any sofic approximation sequence $\Sigma$ of $G$?
\end{question}
\section{Properties of topological pressure}
Let $\alpha$ be a continuous action of a countable sofic group $G$ on a compact metrizable space $X$ and $\Sigma$ a sofic approximation sequence of $G$.
In this section, we study some properties of the map $P_\Sigma(\cdot,X,G): C(X)\rightarrow \Rb\cup \{\pm \infty\}$ and give a sufficient condition involving topological pressure for determining membership in $M_G(X)$ when $G$ is a general countable sofic group.

The following result is well known when $G$ is amenable. For example, see \cite[Theorem 9.7]{Walters82} for the case $G=\Zb$ and \cite[Corollary 5.2.6]{Mou85} for the general case $G$ is amenable.
\begin{proposition}
\label{T-basic properties of sofic topological pressure}
If $f,g\in C(X), s\in G$ and $c\in \Rb$ then the following are true.
\begin{enumerate}
\item [(i)] $P_\Sigma(0,X,G)=h_\Sigma(X,G)$,

\item [(ii)] $P_\Sigma(f+c,X,G)=P_\Sigma(f,X,G)+c$,

\item [(iii)] $P_\Sigma(f+g,X,G)\leq P_\Sigma(f,X,G)+P_\Sigma(g,X,G)$,


\item [(iv)] $f\leq g$ implies $P_\Sigma(f,X,G)\leq P_\Sigma(g,X,G)$. In particular, $h_\Sigma(X,G)+\min f\leq P_\Sigma(f,X,G)\leq h_\Sigma(X,G)+\max f$,

\item [(v)] $P_\Sigma(\cdot, X,G)$ is either finite valued or constantly $\pm\infty$,

\item [(vi)] If $P_\Sigma(\cdot,X,G)\neq\pm\infty,$
then $|P_\Sigma(f,X,G)-P_\Sigma(g,X,G)|\leq \|f-g\|_\infty$, where $\|.\|_\infty$ is the suprenorm on $C(X)$,

\item [(vii)] If $P_\Sigma(\cdot,X,G)\neq\pm\infty$ then $P_\Sigma(\cdot, X,G)$ is convex,

\item [(viii)] $P_\Sigma(f+g\circ\alpha_s-g,X,G)=P_\Sigma(f,X,G)$,

\item [(ix)] $P_\Sigma(cf,X,G)\leq c\cdot P_\Sigma(f,X,G)$ if $c\geq 1$ and $P_\Sigma(cf,X,G)\geq c \cdot P_\Sigma(f,X,G)$ if $c\leq 1$,

\item [(x)] $|P_\Sigma(f,X,G)|\leq  P_\Sigma(|f|,X,G)$.

\end{enumerate}
\begin{proof}
Let $\rho$ be a compatible metric on $X$. Let $\sigma$ be a map from $G$ to $\sym(d)$ for some $d\in\Nb$. Let $\varepsilon,\delta>0$ and $F$ be a nonempty finite subset of $G$.

(i), (ii), (iii) and (iv) are clear from the definition of pressure and Remark \ref{R-topological entropy is a special case of topological pressure}.

(v) From (i) and (ii) we get $P_\Sigma(f,X,G)=\pm\infty$ iff $h_\Sigma(X,G)=\pm\infty$.

(vi) follows from (iii) and (iv).

(vii) By H\"{o}lder's inequality, if $p\in [0,1]$ and $\cE$ is a finite subset of $\map(\rho,F,\delta,\sigma)$ then we have
\begin{eqnarray*}
&&\sum_{\varphi\in\cE} \exp\Big(p\sum_{a=1}^df(\varphi(a))+(1-p)\sum_{a=1}^dg(\varphi(a))\Big)\\
&\leq &\Big(\sum_{\varphi\in\cE} \exp\Big(\sum_{a=1}^df(\varphi(a))\Big)\Big)^p\Big(\sum_{\varphi\in\cE} \exp\Big(\sum_{a=1}^dg(\varphi(a))\Big)\Big)^{1-p}.
\end{eqnarray*}
Therefore, $$M^\varepsilon_{\Sigma,\infty}(pf+(1-p)g,X,G,\rho,F,\delta,\sigma)\leq M^\varepsilon_{\Sigma,\infty}(f,X,G,\rho,F,\delta,\sigma)^p\cdot M^\varepsilon_{\Sigma,\infty}(g,X,G,\rho,F,\delta,\sigma)^{1-p},$$ and (vii) follows.

(viii) Let $\sigma$ be a map from $G$ to $\sym(d)$ for some $d\in\Nb$. Let $\varepsilon,\kappa>0$ and $F$ be a nonempty finite subset of $G$ containing $s$. Since $g$ is continuous there exists $Q>0$ such that $|g(x)|\leq Q$ for any $x\in X$. Choose $\delta>0$ such that $2Q\delta |F|<\kappa$ and $|g(y)-g(z)|<\kappa$ for any $y,z\in X$ with $\rho(y,z)<\sqrt\delta$.
Let $\cE$ be a $(\rho_\infty,\varepsilon)$-separated subset of $\map(\rho,F,\delta,\sigma)$. For each $\varphi\in \cE$ we denote by $\Lambda_\varphi$ the set of all $a\in[d]$ such that $\rho(\varphi(ta),t\varphi(a))<\sqrt{\delta}$ for all $t\in F$. Then $|\Lambda_\varphi|\geq (1-|F|\delta)d$ and so
\begin{eqnarray*}
&&\exp\Big(\sum_{a=1}^d(g(s\varphi(a))-g(\varphi(sa)))\Big)\\
&=&\exp\Big(\sum_{a\in\Lambda_\varphi}(g(s\varphi(a))-g(\varphi(sa)))\Big)\exp\Big(\sum_{a\notin\Lambda_\varphi}(g(s\varphi(a))-g(\varphi(sa)))\Big)\\
&\leq& \exp(\kappa d)\exp(2Q|F|\delta d).
\end{eqnarray*}
Therefore,
\begin{eqnarray*}
&&\sum_{\varphi\in\cE}\exp\Big(\sum_{a=1}^d(f+g\circ\alpha_s-g)(\varphi(a))\Big)\\
&=&
\sum_{\varphi\in\cE}\exp\Big(\sum_{a=1}^df(\varphi(a))\Big)\exp\Big(\sum_{a=1}^d(g(s\varphi(a))-g(\varphi(sa)))\Big)\\
&\leq& \sum_{\varphi\in\cE}\exp\Big(\sum_{a=1}^df(\varphi(a))\Big)\exp(\kappa d)\exp(2Q|F|\delta d).
\end{eqnarray*}
Thus
\begin{eqnarray*}
\log M^\varepsilon_{\Sigma,\infty}(f+g\circ\alpha_s-g,X,G,\rho,F,\delta,\sigma)&\leq& \log M^\varepsilon_{\Sigma,\infty}(f,X,G,\rho,F,\delta,\sigma) +\kappa d+2P|F|\delta d \\
&\leq& \log M^\varepsilon_{\Sigma,\infty}(f,X,G,\rho,F,\delta,\sigma) +2\kappa d,
\end{eqnarray*}
and hence $P^\varepsilon_{\Sigma,\infty}(f+g\circ\alpha_s-g,X,G,\rho,F)\leq P^\varepsilon_{\Sigma,\infty}(f,X,G,\rho,F)+2\kappa$ for any nonempty finite subset $F$ of $G$, $\varepsilon>0$ and $\kappa>0$. Therefore, $P_{\Sigma,\infty}(f+g\circ\alpha_s-g,X,G,\rho)\leq P_{\Sigma,\infty}(f,X,G,\rho)+2\kappa$, for any $\kappa>0$.

Similarly, from
\begin{eqnarray*}
&&\exp\Big(\sum_{a=1}^d(g(s\varphi(a))-g(\varphi(sa)))\Big)\\
&=&\exp\Big(\sum_{a\in\Lambda_\varphi}(g(s\varphi(a))-g(\varphi(sa)))\Big)\exp\Big(\sum_{a\notin\Lambda_\varphi}(g(s\varphi(a))-g(\varphi(sa)))\Big)\\
&\geq& \exp(-\kappa d)\exp(-2Q|F|\delta d),
\end{eqnarray*}
we get $P_{\Sigma,\infty}(f+g\circ\alpha_s-g,X,G,\rho)\geq P_{\Sigma,\infty}(f,X,G,\rho)-2\kappa$, for any $\kappa>0$.
Therefore, $P_{\Sigma,\infty}(f+g\circ\alpha_s-g,X,G,\rho)= P_{\Sigma,\infty}(f,X,G,\rho)$.

(ix) If $a_1,...,a_k$ are positive numbers with $\sum_{i=1}^k a_i=1$ then $\sum_{i=1}^k a_i^c\leq 1$ when $c\geq 1$, and $\sum_{i=1}^k a_i^c\geq 1$ when $c\leq 1$. Hence if $b_1,...,b_k$ are positive numbers then $\sum_{i=1}^kb_i^c\leq \Big(\sum_{i=1}^kb_i\Big)^c$ when $c\geq 1$, and $\sum_{i=1}^kb_i^c\geq \Big(\sum_{i=1}^kb_i\Big)^c$ when $c\leq 1$. Therefore, if $\cE$ is a finite subset of $\map(\rho,F,\delta,\sigma)$ we have $$\sum_{\varphi\in \cE}\exp\Big(c\sum_{j=1}^d f(\varphi(j))\Big)\leq \Big(\sum_{\varphi\in \cE}\exp\Big(\sum_{j=1}^d f(\varphi(j))\Big)\Big)^c \mbox{ when } c\geq 1,$$
and $$\sum_{\varphi\in \cE}\exp\Big(c\sum_{j=1}^d f(\varphi(j))\Big)\geq \Big(\sum_{\varphi\in \cE}\exp\Big(\sum_{j=1}^d f(\varphi(j))\Big)\Big)^c \mbox{ when } c\leq 1,$$
Then (ix) follows.

(x) From (iv) and (ix) we get (x).
\end{proof}
\end{proposition}
Let $\cB(X)$ be the $\sigma$-algebra of Borel subsets of $X$.
Recall that a finite signed measure is a map $\mu:\cB(X)\rightarrow \Rb$ satisfying $$\mu(\bigcup_{i=1}^\infty A_i)=\sum_{i=1}^\infty\mu(A_i),$$
 whenever $\{A_i\}_{i=1}^\infty$ is a pairwise disjoint collection of members of $\cB(X)$.

Now we prove a sufficient condition for a finite signed measure to be a member of $M_G(X)$, using topological pressure. It is known for the case of $\Zb$-actions \cite[Theorem 9.11]{Walters82} and we follow the proof there.
\begin{theorem}
\label{T-topoloigcal pressure determines the members of invariant measures}
Assume that $h_\Sigma(X,G)\neq \pm\infty$. Let $\mu:\cB(X)\rightarrow \Rb$ be a finite signed measure. If $\int_X fd\mu\leq P_\Sigma(f,X,G)$ for all $f\in C(X)$, then $\mu\in M_G(X)$.
\begin{proof}

 Suppose $f\geq 0$. If $\kappa>0$ and $n>0$ we have
\begin{eqnarray*}
\int n(f+\kappa)d\mu&=&-\int -n(f+\kappa)d\mu\geq -P_\Sigma(-n(f+\kappa),X,G)\\
&\geq& -[h_\Sigma(X,G)+\max(-n(f+\kappa))] \mbox{ by Theorem \ref{T-basic properties of sofic topological pressure}(iv)}\\
&=&-h_\Sigma(X,G)+n\min(f+\kappa)\\
&>&0 \mbox{ for large n.}
\end{eqnarray*}
Therefore $\int (f+\kappa)d\mu>0$ and hence $\int f d\mu\geq 0$. Thus $\mu$ takes only non-negative values.

If $n\in \Zb$ then $\int nd\mu\leq P_\Sigma(n,X,G)=h_\Sigma(X,G)+n$, so that $\mu(X)\leq 1+h_\Sigma(X,G)/n$ if $n>0$ and hence $\mu(X)\leq 1$, and $\mu(X)\geq 1+h_\Sigma(X,G)/n$ if $n<0$ and hence $\mu(X)\geq 1$. Therefore $\mu(X)=1$.

Lastly we show $\mu\in M_G(X)$. Let $s\in G,n\in\Zb$ and $f\in C(X)$. By Proposition \ref{T-basic properties of sofic topological pressure} (viii), one has $n\int (f\circ \alpha_s -f)d\mu\leq P_\Sigma(n(f\circ \alpha_s -f),X,G)=h_\Sigma(X,G)$. If $n>0$ then dividing both sides by $n$ and letting $n$ go to $\infty$ yields $\int (f\circ \alpha_s -f)d\mu\leq 0$, and if $n<0$ then dividing both sides by $n$ and letting $n$ go to $-\infty$ yields $\int (f\circ \alpha_s -f)d\mu\geq 0$. Therefore $\int f\circ \alpha_sd\mu=\int fd\mu,$ for any $f\in C(X),s\in G$. Thus $\mu\in M_G(X)$.
\end{proof}
\end{theorem}
In the case $G$ is amenable, as a consequence of the variational principle for topological pressure, the converse of Theorem \ref{T-topoloigcal pressure determines the members of invariant measures} is also true; see for example \cite[Theorem 9.11]{Walters82} for the case $G=\Zb$. Thus, it is natural to ask the following question
\begin{question}
Let a countable sofic group $G$ act continuously on a compact metrizable space $X$, $\Sigma$ a sofic approximation sequence of $G$ and $\mu\in M_G(X)$. Do we have $$\int_X fd\mu\leq P_\Sigma(f,X,G), \mbox{ for all } f\in C(X)?$$
\end{question}
 Indeed, when $G$ is a general countable sofic group, we only need to consider the case $h_{\Sigma,\mu}(X,G)=-\infty$ since if $h_{\Sigma,\mu}(X,G)\neq -\infty$ then by Theorem \ref{T-variational principle for topological pressure} we obtain $\int_X fd\mu\leq P_\Sigma(f,X,G),$ for all $f\in C(X)$.

\end{document}